\documentclass{article}

\usepackage{arxiv}

\usepackage[utf8]{inputenc} 
\usepackage[T1]{fontenc}    
\usepackage{hyperref}       
\usepackage{url}            
\usepackage{booktabs}       
\usepackage{amsfonts}       
\usepackage{nicefrac}       
\usepackage{microtype}      
\usepackage{doi}
\usepackage{enumitem}
\usepackage{style}
\usepackage{graphicx}
\usepackage{tikz}
\usepackage{doi}
\usepackage[normalem]{ulem}

\title{The Zipped Finite Element Method: High-order Shape Functions for Polygons}

\author{ \href{https://orcid.org/0000-0001-8642-4258}{\includegraphics[scale=0.06]{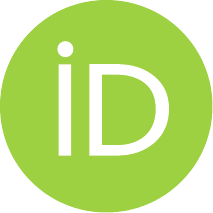}\hspace{1mm}Stefano~Berrone} \\
	Dipartimento di Scienze Matematiche\\
	``G. L. Lagrange''\\
	Politecnico di Torino, TO, 10129 \\
	\texttt{stefano.berrone@polito.it} \\
	\And
	\href{https://orcid.org/0009-0009-6203-9596}{\includegraphics[scale=0.06]{orcid.pdf}\hspace{1mm}Lorenzo~Neva} \\
	Dipartimento di Scienze Matematiche\\
	``G. L. Lagrange''\\
	Politecnico di Torino, TO, 10129 \\
	\texttt{lorenzo.neva@polito.it} \\
	\And
	\href{https://orcid.org/0000-0002-2837-2792}{\includegraphics[scale=0.06]{orcid.pdf}\hspace{1mm}Moreno~Pintore} \\
	Sorbonne Universit\'e, \\
	MEGAVOLT Team, Inria, \\
	4 place Jussieu, 75005 Paris, France \\
	\texttt{moreno.pintore@sorbonne-universite.fr}  \\
	\And
	\href{https://orcid.org/0000-0002-8540-3639}{\includegraphics[scale=0.06]{orcid.pdf}\hspace{1mm}Gioana~Teora} \\
	Dipartimento di Scienze Matematiche\\
	``G. L. Lagrange''\\
	Politecnico di Torino, TO, 10129 \\
	\texttt{gioana.teora@polito.it} \\
	\And
	\href{https://orcid.org/0000-0001-7123-9199}{\includegraphics[scale=0.06]{orcid.pdf}\hspace{1mm}Fabio~Vicini} \\
	Dipartimento di Scienze Matematiche\\
	``G. L. Lagrange''\\
	Politecnico di Torino, TO, 10129 \\
	\texttt{fabio.vicini@polito.it} \\
}

\renewcommand{\shorttitle}{Zipped Finite Element Method}
\hypersetup{
pdftitle={},
pdfsubject={},
pdfauthor={Stefano~Berrone, Lorenzo~Neva, Moreno~Pintore, Gioana~Teora, Fabio~Vicini},
pdfkeywords={},
}

\begin{document}
\maketitle

\begin{abstract}
In this paper, we present a new polygonal finite element method, called the Zipped Finite Element Method, for star-shaped polygons. The proposed approach constructs high-order shape functions as linear combinations of standard finite element basis functions defined on a local trivial sub-triangulation of each element. This refinement is used solely for the construction of the shape functions and does not affect the final number of degrees of freedom. The resulting finite element space includes polynomials of the desired order and preserves conformity across elements. Consequently, the method inherits the convergence properties of the finite element framework under suitable mesh assumptions. Numerical experiments confirm the expected rates of convergence.
\end{abstract}


\section{Introduction}

The simulation of physical phenomena plays a crucial role in many engineering applications, where the numerical solution of partial differential equations is often required. Among the available numerical techniques, the Finite Element Method (FEM) \cite{ciarlet2002finite, brennerscott} is one of the most widely used, due to its flexibility and solid theoretical foundations. However, standard FEM formulations are constrained by the requirements of domain discretization. In two-dimensional problems, such as those considered in this work, the mesh is typically restricted to triangular or convex quadrilateral elements. 
However, in many practical applications, the ability to include polygons with an arbitrary number of vertices, concave elements, or elements with aligned edges can be highly advantageous.
Such flexibility, for instance, can be beneficial to accurately discretize complex geometries \cite{GrappeinTeora2025, Benedetto2017} or to enable efficient refinement and coarsening strategies \cite{Vicini2024}.

Over the years, several approaches have been proposed to extend the FEM to more general meshes. Among the most prominent is the Virtual Element Method (VEM) \cite{LBe13, LBe14}, which uses projection and stabilization operators to define a computable bilinear form without requiring explicit expressions for the basis functions.
Other recent methods have focused on deriving closed-form expressions for local basis functions, often as approximations of VEM functions. Examples include the Lightning Virtual Element Method (LVEM) \cite{TrezziZerbinati2024}, which leverages the Laplace solver from \cite{Gopal2019}, and the Neural Approximated Virtual Element Method (NAVEM) \cite{PintoreTeora2025, navemElasticity}, which employs artificial neural networks.

Despite the evident complexity of defining high-order basis functions over generic polygons, their superior accuracy and convergence properties have long motivated research efforts in this direction.
Building on this idea, the objective of this work is to design an algorithm that constructs high-order shape functions for meshes composed of general star-shaped elements. The proposed method handles convex and concave elements, as well as elements with hanging nodes, using a unified numerical treatment.

The method introduced in this paper is called the Zipped Finite Element Method (Z-FEM in short). The term zipped reflects the key design principle of the approach: local discrete spaces are ``compressed'' to contain a number of degrees of freedom that corresponds with the serendipity ``lazy'' choice \cite{DaVeigaBrezzi2016} and that is sufficient to ensure conformity and preserve polynomial consistency over generic polygons. 

The key idea of Z-FEM is to define high-order local basis functions as weighted combinations of standard FEM shape functions defined on a simple sub-triangulation of each element. This sub-triangulation, usually required for numerical integration, is obtained by connecting the vertices of the polygon to a suitably chosen central point located within the kernel of the element. The coefficients of the linear combinations are determined by solving a local optimization problem, which can be efficiently decomposed into smaller sub-problems thanks to its structure. Moreover, the local space is designed so that polynomials of degree up to the method order can be reproduced. This property is fundamental, as it allows the method to inherit the theoretical framework of the FEM, ensuring both well-posedness of the discrete problem and optimal a priori error estimates.

To the best of our knowledge, several related works have explored similar ideas to those presented in this paper, but remain limited to low-order formulations. In \cite{Bishop2014}, the authors introduce a method based on conforming shape functions constructed as weighted combinations of standard FEM shape functions defined over a fan of simplices forming the polyhedron. Their formulation, restricted to first-order accuracy, is applied to star-shaped polygons in nonlinear solid mechanics. A first-order formulation is also proposed in \cite{Bunge2020} and later extended to second order in \cite{Bunge2022}, where the authors emphasize the challenges of generalizing such basis functions to higher-order cases. Higher-order extensions of the FEM have also been proposed, e.g., in \cite{arbogast2022direct, arbogast2023direct}, but they remain limited to convex elements.

The paper is structured as follows. Section~\ref{sec:notationsmodel} presents the model problem and the notation used in the manuscript. Section~\ref{sec:zippedfem} details and analyzes the required properties of the shape functions. In Section~\ref{sec:shapefunctions} we provide a theoretical justification for the appropriate choice of degrees of freedom. The optimization problem is presented, along with the algorithm adopted for its resolution. In Section~\ref{sec:discreteproblem} we define the local discrete space, demonstrate its local unisolvence, and discuss the theoretical properties of the new proposed method. Section~\ref{sec:numericalexperiments} is dedicated to numerical experiments, which show the effectiveness of the method. Conclusions can be found in Section~\ref{sec:conclusion}.

\section{Notations and the Model problem}
\label{sec:notationsmodel}

Throughout the paper, the usual notation for Sobolev spaces is adopted. Let $\genericset \subset \R^2$ be a bounded open domain. Given two scalar functions $f,g \in L^2(\genericset)$ and two vector fields $\mathbf{s}, \mathbf{t} \in \big[ L^2(\genericset)\big]^2$, we denote by $\scal[\genericset]{f}{g}$ and $\scal[\genericset]{\mathbf{s}}{\mathbf{t}}$ the two bilinear forms
\begin{equation*}
    \scal[\genericset]{f}{g} := \int_\genericset fg \operatorname{d} \omega, \quad \quad \scal[\genericset]{\mathbf{s}}{\mathbf{t}} := \int_\genericset \mathbf{s} \cdot \mathbf{t} \operatorname{d} \omega.
\end{equation*}
Furthermore, $\norm[\sob{s}{\genericset}]{\cdot}$ and $\seminorm[\sob{s}{\genericset}]{\cdot}$ denote the norm and seminorm of functions in $\sob{s}{\genericset}$, for $s \geq 0$, respectively. By definition, $\sob{0}{\genericset} \equiv \leb{2}{\genericset}$.    Lastly, the notation $\norm[2]{}$ denotes the standard Euclidean norm on $\R^2$.

Given a generic polygon $E$, let us introduce the space $\Poly{k}{E}$ of two-dimensional polynomials of degree up to $k \in \N$ on $E$, whose dimension is $n_k := \frac{(k+1)(k+2)}{2}$. We adopt the standard convection $\Poly{-m}{E} = \emptyset$ and $n_{-m} = 0$ for all $m \in \N$.

Let $\Omega \subset \R^2$ be a bounded domain with Lipschitz boundary $\partial \Omega$. We consider the following scalar diffusion-reaction problem with homogeneous Dirichlet boundary conditions:
\begin{equation}
    \label{eq:modelproblem}
    \begin{cases}
        -\nabla \cdot (\D \nabla u) + \gamma u = f & \text{in } \Omega,\\
        u = 0 &\text{on } \Gamma = \partial \Omega,\\
    \end{cases}
\end{equation}
where $\D \in [\leb{\infty}{\Omega}]^{2 \times 2}$ is a diffusion tensor, uniformly symmetric positive definite over $\Omega$, i.e. there exist two constants $0 < \alpha_1 \leq \alpha_2$ such that 
\begin{equation}
    \alpha_1 \norm[2]{\bm{\epsilon}} \leq \bm{\epsilon}^T \D(\xx) \bm{\epsilon} \leq \alpha_2 \norm[2]{\bm{\epsilon}} \qquad \forall \xx \in \Omega,\quad \forall\bm{\epsilon} \in \R^2.
\end{equation}
Moreover, $\gamma \in \leb{\infty}{\Omega}$, $\gamma(\xx) \geq 0$ for all $\xx \in \Omega$, is a reaction coefficient, and $f \in \leb{2}{\Omega}$ is a scalar loading term.

The variational formulation of Problem \eqref{eq:modelproblem} reads as: \textit{Find $u \in \VP := \sob[0]{1}{\Omega}$ such that}
\begin{equation}
    \mathcal{B}(u, v) = \scal[\Omega]{f}{v} \quad \forall v \in \VP,
    \label{eq:varproblem}
\end{equation}
where $\bilin{}{} : \VP \times \VP \to \R$ is the following symmetric, continuous, and coercive bilinear form
\begin{equation}
\label{eq:formabilincontinua}
\mathcal{B}(u,v) = \scal[\Omega]{\D \nabla u}{\nabla v} + \scal[\Omega]{\gamma u}{v} \quad \forall u,v \in \VP.
\end{equation}
Thanks to the Lax-Milgram Theorem \cite{Evans2010}, the Problem \eqref{eq:varproblem} admits a unique solution.

\section{The Zipped Finite Element Space}
\label{sec:zippedfem}

Let $\Th$ be a discretization of the domain $\Omega$ into non-overlapping polygonal elements $E$. The symbols $\Nv[E]$, $\Eh[E]$, and $h_E$ denote the number of vertices, the set of edges, and the diameter of the element $E$, respectively. Moreover we set $h := \max_{E \in \Th} h_E$.

We assume $\Th$ to satisfy the following mesh assumptions.
\begin{assumptions}[Mesh assumptions]
\label{ass:meshassumption}
There exists a positive real number $\rho \in (0, 1)$, independent of $h$, such that
    \begin{itemize}
        \item for each element $E\in\Th$ and each edge $e \in \Eh[E]$, it holds $\vert e \vert \geq \rho h_E$;
        \item each element $E\in\Th$ is star-shaped with respect to a ball of radius $\geq \rho h_E$.
    \end{itemize}
\end{assumptions}

Given an integer number $k \geq 1$, we define a finite-dimensional space $\VPh{h,k} \subset \VP$ as
\begin{equation}
    \VPh{h,k} := \{v \in \con{0}{\overline{\Omega}} \cap \VP:\ v_{|E} \in \VPh[E]{k} \ \forall E \in \Th\},
    \label{eq:global_space}
\end{equation}
where the elemental space $\VPh[E]{k}$ can be described as
\begin{equation}
    \VPh[E]{k} := \myspan\{\varphi_i: i = 1,\dots, \Ndof[E]\} \quad \text{with}\quad \Ndof[E] = \dim \VPh[E]{k},
    \label{eq:ele_space}
\end{equation}
for a given set of basis functions $\{\varphi_i\}_{i=1}^{\Ndof[E]}$. We aim to define the shapes of these basis functions and their number $\Ndof[E]$ so that the following properties hold:
\begin{enumerate}[label=\textbf{P.\arabic*}]
    \item \label{prop:polynomial}The set $\VPh[E]{k}$ contains the set of polynomials $\Poly{k}{E}$. 
    \item \label{prop:krnocker}The elemental shape functions satisfy the Kronecker-Delta property with respect to a set of linearly independent linear operators
    \begin{equation}
        \dof_i : \VPh[E]{k} \to \R\qquad i = 1,\dots,\Ndof[E],
        \label{eq:lineardofsoperator}
    \end{equation}
    termed Degrees of Freedom (DOFs in short), i.e.
    $\dof_i(\varphi_j) = \delta_{ij}$ for all $i,j=1,\dots,\Ndof[E]$. 
    \item\label{prop:continuity} The basis functions are continuous across adjacent elements.
\end{enumerate}

For this purpose, the first issue concerns the choice of the number of degrees of freedom that enable us to satisfy Properties \ref{prop:polynomial} and \ref{prop:continuity}. This problem has been extensively investigated in the framework of serendipity element methods \cite{arbogast2023direct, DaVeigaBrezzi2016}.
The number of degrees of freedom $\Ndof[E]$ must be at least $n_k$, to include all polynomials of order $\leq k$ and satisfy Property \ref{prop:polynomial}.
Moreover, to guarantee Property \ref{prop:continuity}, $(k + 1)$ degrees of freedom must be placed along each edge, resulting in a total number of $\Nv[E] k $ degrees of freedom on the boundary $\partial E$ of the element $E$. In particular, if $k \leq 2$, the boundary degrees of freedom are sufficient to ensure both conformity and the unique identification of polynomials, since $\Nv[E] \geq 3$. When $k \geq 3$, instead, the minimum number of degrees of freedom required to satisfy Properties \ref{prop:polynomial} and \ref{prop:continuity} for a convex polygon $E$ is given by \cite{arbogast2023direct, DaVeigaBrezzi2016}
\begin{equation*}
    \Nv[E]k + n_{(k - \eta_{E})},
\end{equation*}
where $\eta_E$ denotes the number of straight lines needed to cover the entire boundary $\partial E$. For stability reasons, however, it is also important to determine the number $\eta_E$ such that consecutive edges lying \textit{almost} on the same straight lines are counted as one, while carefully defining what ``almost'' means. 
To avoid the stability and geometric issues arising from this ambiguity, while being able to handle convex, concave, and hanging-node elements indistinctly, we adopt the so-called \emph{serendipity lazy choice} \cite{DaVeigaBrezzi2016}, which consists of always including $n_{k-3}$ internal degrees of freedom, that are sufficient to guarantee both stability and geometric independence. 

Thus, the final number of local degrees of freedom is
\begin{equation}
    \Ndof[E]  = \Nv[E] k + n_{k-3}.
    \label{eq:totalnumberdofs}
\end{equation}


\section{High-order shape functions on polygons}\label{sec:shapefunctions}

Let us consider a star-shaped polygon $E \subset \R^2$. Our aim is to determine a set of $\Ndof[E]$ shape functions $\{\varphi_i\}_{i=1}^{\Ndof[E]}$ on $E$ for a generic order $k \geq 1$ such that they satisfy Properties \ref{prop:polynomial}, \ref{prop:krnocker}, and\ref{prop:continuity}. 

For this purpose, we adopt a similar approach to that proposed by \cite{Bunge2020} for the lowest order $k=1$, while extending it to a generic order $k \geq 1$. Specifically, we compute these shape functions as a weighted combination of Finite Element shape functions living on a sub-triangulation of the element $E$, which can be trivially built under the assumption that $E$ is star-shaped. 
Indeed, given an internal point $\xx^E$ of $E$, we can build a sub-triangulation $\TE[E] = \{T_j\}_{j=1}^{\Nv[E]}$ of $E$ by connecting the $\Nv[E]$ vertices $\{\vv_j\}_{j=1}^{\Nv[E]}$ of $E$ with $\xx^E$ to form the following triangles: 
\begin{equation}
    T_j = \{\xx^E, \vv_j,\vv_{j+1}\}, \quad  j=1,\dots, \Nv[E],
    \label{eq:triangles}
\end{equation}
where the symbol $\vv_{\Nv[E]+1}$ denotes the vertex $\vv_1$.

In particular, we define the point $\xx_E$ as the center of the disk of radius $r_E$ with respect to which the polygon $E$ is star-shaped. This point can be found as the solution of the following linear programming problem \cite{Calafiore2014}:
\begin{equation}
\begin{aligned}
    \displaystyle\max_{(\xx^E, r^E) \in \R^3}\quad &r^E\\
    \text{such that} \quad& \nn_e \cdot \xx^E + r^E \leq b_e \quad \forall e \in \Eh[E]\\
    & 0 < r_E \leq h_E.
\end{aligned}
\label{eq:opt:centerstarshaped}
\end{equation}
Moreover, for each edge $e \in \Eh[E]$, $\nn_e$ denotes the unit outward normal vector to the edge $e$, and the scalar coefficient $b_e$ is chosen in such a way that the edge lies on the line defined by the equation $\nn_e \cdot \xx = b_e$ for all $\xx \in e$.

\begin{remark}\label{rem:shaperegular}
Another possible choice for $\xx^E$ is the one described in \cite{Bunge2020}, which yields highly accurate results in the linear case. Specifically, in such a paper, the internal point is determined as a weighted combination of the element vertices, i.e. $\xx^E := \sum_{j=1}^{\Nv[E]} a_j \vv_j$, where the weights are the least-norm solution of the following optimization problem: 
\begin{equation}
    \begin{aligned}
        \min_{\{a_j\}_{j=1}^{\Nv[E]}} \quad&\sum_{j=1}^{\Nv[E]} \mathrm{area}(T_j)^2\\
        \text{such that} \quad&\sum_{j=1}^{\Nv[E]} a_j = 1.
    \end{aligned}
    \label{eq:opt:ssa}
\end{equation}
Here, instead, we adopt the alternative choice determined by \eqref{eq:opt:centerstarshaped} because it guarantees that the resulting global triangulation $\bigcup_{E \in \Th}\TE[E]$ is shape-regular \cite{Brenner2017}. This property allows us to establish optimal error estimates (see Section \ref{sec:discreteproblem}). Nonetheless, numerical experiments do not reveal any significant difference in accuracy when using the various choices for $\xx^E$. In Figure \ref{fig:differenttriangul}, we illustrate the sub-triangulations obtained by defining $\xx^E$ as the solution of \eqref{eq:opt:ssa} and \eqref{eq:opt:centerstarshaped} for a heptagon.

\begin{figure}[!ht]
    \centering
    \begin{subfigure}{0.45\textwidth}
        \includegraphics[width=1\linewidth]{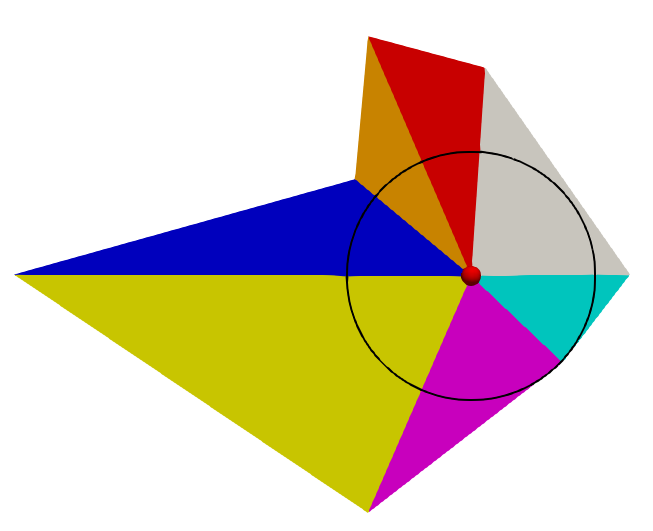}
        \caption{}
    \end{subfigure}
    \begin{subfigure}{0.45\textwidth}
        \includegraphics[width=1\linewidth]{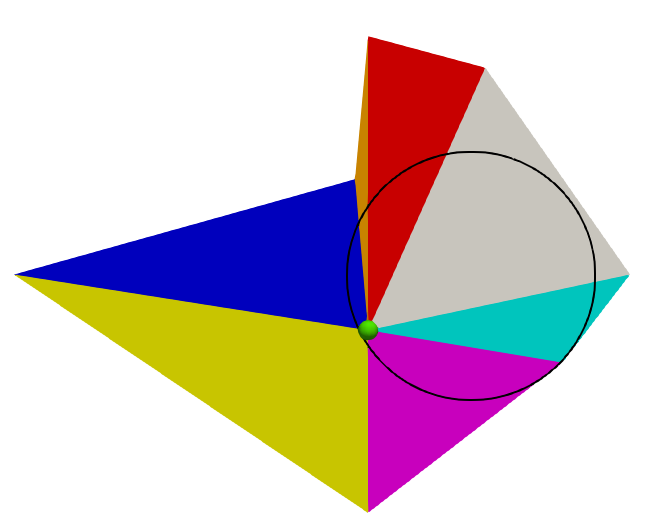}
        \caption{}
    \end{subfigure}
    \caption{Different sub-triangulations of a heptagon. Left: Sub-triangulation obtained by defining $\xx^E$ (the red dot) as the solution of \eqref{eq:opt:centerstarshaped}. Right: Sub-triangulation obtained by defining $\xx^E$ (the green dot) as the solution of \eqref{eq:opt:ssa}. The black line denotes the circumference of the ball with respect to which the heptagon is star-shaped.}
    \label{fig:differenttriangul}
\end{figure}
\end{remark}

\subsection{Standard basis functions over the element sub-triangulation}

The standard finite element space of order $k \geq 1$ over the sub-triangulation $\TE[E]$ is defined as
\begin{equation}
    \mathbb{V}_k(E; \TE[E]) := \{v \in \con{0}{\overline{E}} \cap \sob{1}{E} : v_{|T} \in \Poly{k}{T}\ \ \forall\, T \in \TE[E] \}.
    \label{eq:fintielementlocalspace}
\end{equation}
We observe that $\Poly{k}{E} \subset \mathbb{V}_k(E; \TE[E])$. Moreover, the dimension of this space is given by
\begin{equation*}
    N^E_{\Psi,k} = \left[(2k - 1)  +  n_{k-3} \right] \Nv[E]+ 1.
\end{equation*} 

Let us now introduce the finite element shape functions $\{\widehat{\Psi}_{\ell}\}_{\ell = 1}^{n_k}$ over the reference triangle $\widehat{T} = \{(0,0), (1,0), (0,1)\}$. 
Given the vector-valued function $\mathbf{pos} : \{1,\dots,n_k\} \to \N^2$ defined as 
\begin{equation*}
    \mathbf{pos}(1) = (0, 0),\quad \dots,\quad \mathbf{pos}(k) = (k, 0),\quad \mathbf{pos}(k + 1) = (0, 1),\quad \dots, \quad \mathbf{pos}(2k - 1) = (k-1, k-1),\ \dots,
\end{equation*}
we introduce the following set of nodal coordinates in $\widehat{T}$:
\begin{equation*}
    \widehat{\bm{\xi}}_{\ell} = \frac{1}{k} \mathbf{pos}(\ell), \quad \text{for } \ell = 1, \dots,n_k, 
\end{equation*}
which are evenly distributed in each dimension. 
The finite element shape functions $\{\widehat{\Psi}_{\ell}\}_{\ell = 1}^{n_k}$ defined with respect to this set of nodes are given by
\begin{equation*}
    \widehat{\Psi}_{\ell}(\widehat{\xx}) = \prod_{d=1}^3 \prod_{m = 0}^{\mathbf{pos}(\ell) - 1} \frac{\lambda_d(\widehat{\xx}) - \frac{m}{k}}{\lambda_d(\widehat{\bm{\xi}}_{\ell})- \frac{m}{k}} \qquad \forall\, \widehat{\xx} := (\widehat{x}_1, \widehat{x}_2) \in \R^2,
\end{equation*}
where $\lambda_d$ represent the barycentric coordinates related to this reference elements, i.e.
\begin{equation*}
    \lambda_d(\widehat{\xx}) = \widehat{x}_d\quad d=1,2\quad\text{and}\quad \lambda_3(\widehat{\xx}) = 1 -\sum_{d=1}^2 \widehat{x}_d.
\end{equation*}
We properly map these finite element functions to the sub-triangles $\{T_j\}_{j=1}^{\Nv[E]}$, and we glue them together to obtain a set of basis functions $\{\Psi_n\}_{n=1}^{N^E_{\Psi,k}}$ for $\mathbb{V}_k(E; \TE[E])$ that satisfy the Kronecker-Delta property with respect to the full set of nodes $\{\bm{\xi}_n\}_{n=1}^{N^E_{\Psi,k}}$, i.e.
\begin{equation}
    \label{eq:nodes}
    \Psi_m(\bm{\xi}_n) = \delta_{nm} \quad \forall n,m \in \finenodes[E] := \{ 1,\dots,N^E_{\Psi,k} \}.
\end{equation}

\subsection{Selection of the nodal coordinates}

We observe that this set of functions satisfies properties~\ref{prop:polynomial}, \ref{prop:krnocker}, and \ref{prop:continuity}. Nevertheless, to reduce the dimension of the final space $\VPh[E]{k}$, we select a subset $\coarsenodes[E] \subset \finenodes[E]$ of $\Ndof[E]$ nodes $\{\xx_i\}_{i \in \coarsenodes[E]}$ that allows us to satisfy these properties, while satisfying
\begin{equation}
    \Ndof[E] \ll N^E_{\Psi,k}.
\end{equation}
Here, $\Ndof[E]$ is defined as in Equation~\eqref{eq:totalnumberdofs}. Moreover, let us denote by $\virtualnodes[E] := \finenodes[E] \setminus \coarsenodes[E]$ the set of remaining indices, associated with the nodes $\{\pp_j\}_{j \in \virtualnodes[E]} \subset \{\bm{\xi}_i\}_{i \in \finenodes[E]}$, whose cardinality is
\begin{equation*}
    \Nvn[E] = N^E_{\Psi,k} - \Ndof[E].
\end{equation*}
In the following, we will name the sets $\{\bm{\xi}_{n}\}_{n \in \finenodes[E]}$, $\{\xx_{i}\}_{i \in \coarsenodes[E]}$, and $\{\pp_{j}\}_{j \in \virtualnodes[E]}$ as the set of fine, coarse, and virtual nodes, respectively.
Moreover, without loss of generality, we always assume that $\coarsenodes[E] = \{1,\dots,\Ndof[E]\}$ and, as a consequence, $\virtualnodes[E] = \{\Ndof[E] + 1, \dots, N^E_{\Psi,k}\}$ to simplify notations.

Given this set of coarse nodes $\{\xx_i\}_{i \in \coarsenodes[E]}$, we define the degrees of freedom introduced in Equation~\eqref{eq:lineardofsoperator} as
\begin{equation}
    \label{eq:dofs}
    \dof_i(v) = v(\xx_i) \quad \forall i \in \coarsenodes[E] \quad \text{ and }\quad \forall v \in \VPh[E]{k}.
\end{equation}

The set of coarse nodes with respect to which we define these linear operators must include all $\Nv[E] k$ degrees of freedom located on the boundary $\partial E$ to ensure Property~\ref{prop:continuity}. In addition, $n_{k-3}$ DOFs internal to $E$ must also be included (see the definition of $\Ndof[E]$ in \eqref{eq:totalnumberdofs}). There is no prescribed choice for these internal degrees of freedom, but a strategy to properly select such internal nodal coordinates to guarantee polynomial reproducibility may be based on the following well-established theorem. 

\begin{theorem}[Theorem 6.1 in \cite{Gregory2007}]\label{theor:ChungYao}
    Suppose $\{L_0,\dots,L_k\}$ are $k+1$ distinct lines on $\R^2$ and $U = \{\xx_1, \dots, \xx_{n_k}\}$ is a set of $n_k$ distinct points such that $\xx_1 \in L_0$, $\xx_2, \xx_3 \in L_1\setminus L_0$, $\dots$, and $\xx_{n_k - k},\dots,\xx_{n_k} \in L_k \setminus \{L_0,\dots,L_{k-1}\}$. Then, there exists a unique polynomial of degree $\leq k$ that interpolates on $U$. Moreover, if $X = \{\xx_1,\dots,\xx_K\}$ with $U \subset X$, then $X$ is $k$-unisolvent, i.e., the values of any polynomial of degree $\leq k$ at the points of $X$ uniquely determine the polynomial itself.
\end{theorem}

From this theorem, one can establish a rigorous criterion for selecting at least $n_{k-3}$ internal points among the available degrees of freedom to ensure the inclusion of the required polynomials. 

\begin{figure}[!ht]
    \centering
    \begin{subfigure}{0.25\textwidth}
        \includegraphics[width=1\linewidth]{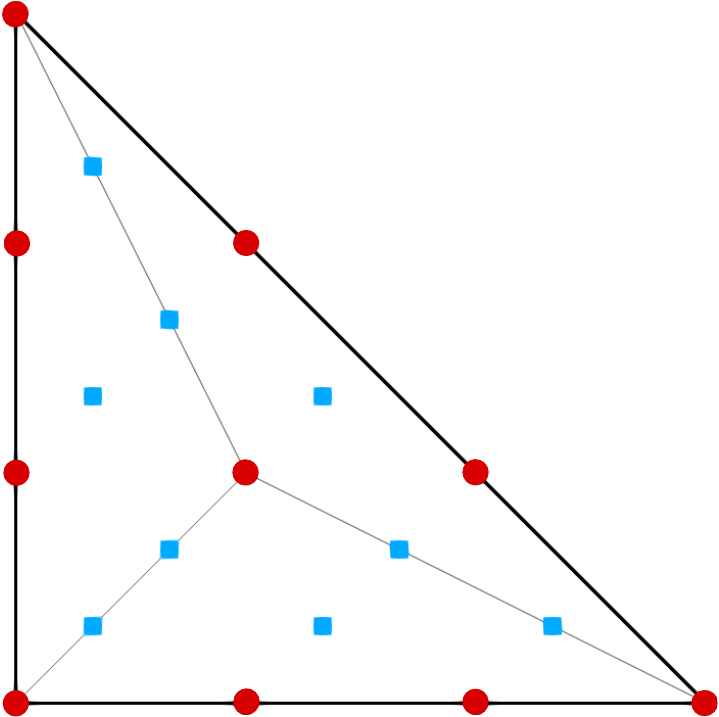}
        \caption{$k = 3$.}
        \label{fig:fem_dofs:3}
    \end{subfigure}\hspace{10pt}
    \begin{subfigure}{0.35\textwidth}
        \includegraphics[width=1\linewidth]{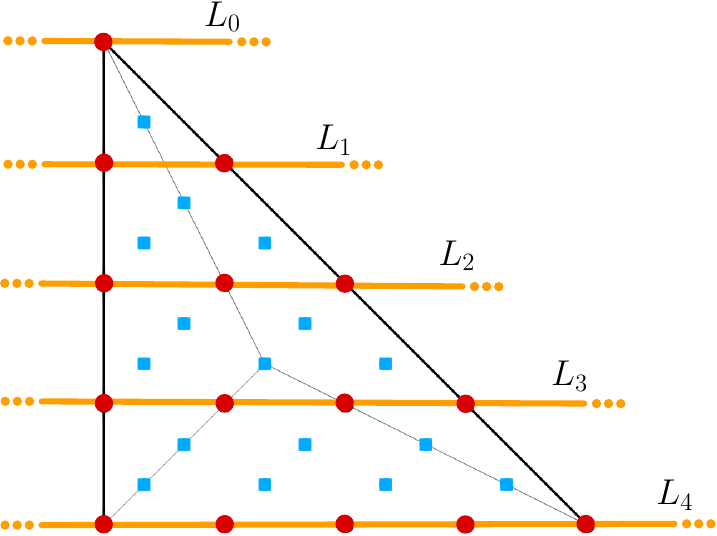}
        \caption{$k = 4$.}
        \label{fig:fem_dofs:4}
    \end{subfigure}\hspace{10pt}
    \begin{subfigure}{0.25\textwidth}
        \includegraphics[width=1\linewidth]{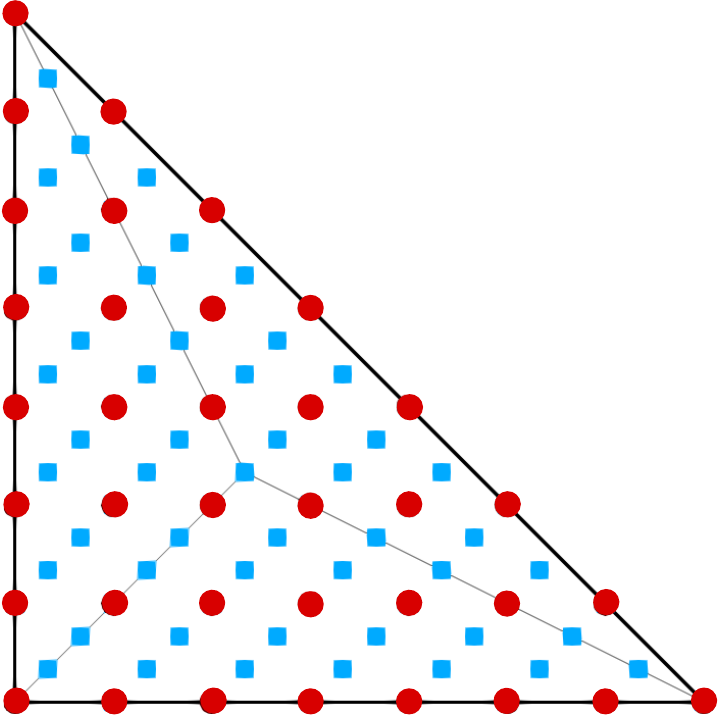}
        \caption{$k = 7$.}
        \label{fig:fem_dofs:7}
    \end{subfigure}
    \caption{Selection of local $\Ndof[E]$ DOFs on a triangular element $E$ for different orders $k$. Red dots indicate $\coarsenodes[E]$, while blue squares represent $\virtualnodes[E]$. The yellow lines on $k=4$ represent the distinct lines $\{L_0, \dots, L_k\}$ defined in Theorem~\eqref{theor:ChungYao}.}
    \label{fig:fem_dofs}
\end{figure}

In particular, we observe that the standard Finite Element Method on triangular meshes can be derived from our method using the rule provided by Theorem~\ref{theor:ChungYao} for the selection of the coarse nodes. To illustrate this, Figure~\ref{fig:fem_dofs} shows the set of fine nodes on the reference triangle for orders $k \in \{3, 4, 7\}$, highlighting the choice for the coarse nodes for which the Z-FEM space coincides exactly with the standard Lagrange FEM space. 
For $k > 2$, the selected points are located along the equispaced lines $\{L_0, \dots, L_k\}$ parallel to one of the edges of the $E$, which are drawn for the case $k=4$ in Figure~\ref{fig:fem_dofs:4}. We remark that this choice satisfies Theorem~\ref{theor:ChungYao}, while allowing the Z-FEM formulation to coincide exactly with the classical Lagrange finite element on triangular meshes. 

However, when more intricate geometries are handled, such as concave elements or elements having hanging nodes, selecting internal points that rigorously satisfy Theorem~\ref{theor:ChungYao} becomes geometrically complex and computationally demanding. 
For these reasons, we propose an alternative strategy for the internal selection of points aimed at reducing the local computational cost.  
In the numerical tests, this strategy is applied to all polygons independently of their geometry and leads to satisfactory results in each tested case.  

The proposed selection strategy is guided by two main principles:
\begin{itemize}
    \item minimizing internal point alignment;
    \item ensuring a homogeneous spatial distribution of internal points within the polygon.
\end{itemize}

\begin{figure}[!ht]
    \centering
    \begin{subfigure}{0.4\textwidth}
        \includegraphics[width=1\linewidth]{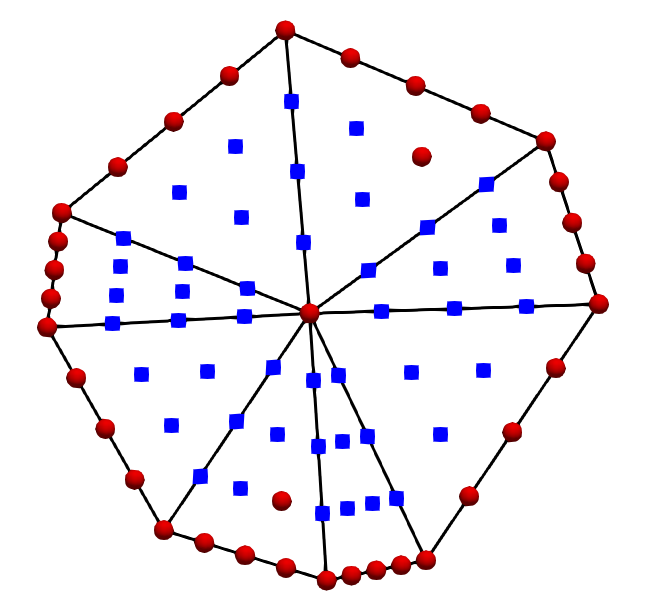}
        \caption{}
    \end{subfigure}
    \begin{subfigure}{0.4\textwidth}
        \includegraphics[width=1\linewidth]{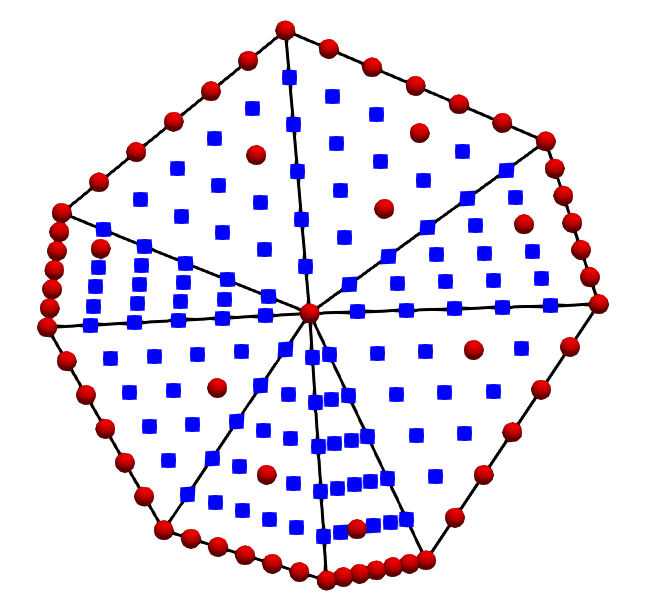}
        \caption{}
    \end{subfigure}
    \caption{Selection of local $\Ndof[E]$ DOFs with the proposed heuristic procedure on a convex octagon $E$. Red dots indicate $\coarsenodes[E]$, while blue squares represent $\virtualnodes[E]$. Left: $k=4$. Right: $k=6$.}
    \label{fig:poly_dofs}
\end{figure}

To this end, a heuristic reordering algorithm is implemented. For any $k \geq 3$, for simplicity, we always select the nodal coordinates corresponding to $\xx_E$ among the coarse nodes. Thus, for $k \geq 4$, we have to select the remaining $n_{k-3} - 1$ internal coarse nodes. Starting from the natural vertex ordering, the triangles are divided into $n_{k-3} - 1$ groups and then interleaved following a regular pattern determined by the ratio between the total number of triangles and the number of groups.
This procedure yields a spatially balanced indexing that provides uniform coverage of the domain and prevents local clustering of points.
The same reordering pattern is then applied to determine the sequence of internal point selection within each triangle, further reducing local alignment and improving geometric uniformity. The subdivision of the triangles into groups and the corresponding assignment of the DOFs to each triangle in the case of an octagon is shown in the Table~\ref{tab:dofassignment} for $k \in \{4,5,6,7\}$, while the corresponding selected nodal coordinates are shown in the Figure~\ref{fig:poly_dofs} for $k \in \{4,6\}$.

\begin{table}[!ht]
\centering
\caption{Number of internal degrees of freedom selecting for each triangle inside an octagon according to our heuristic strategy. Vertical lines denote triangles group subdivision.}
\label{tab:dofassignment}
\begin{tabular}{cccccccccc}
$k$ & $n_{k-3}-1$             & $T_1$                  & $T_2$                  & $T_3$                  & $T_4$                  & $T_5$                  & $T_6$                  & $T_7$                  & $T_8$ \\ \hline
4   & \multicolumn{1}{c|}{2}  & 1                      & 0                      & 0                      & \multicolumn{1}{c|}{0} & 1                      & 0                      & 0                      & 0     \\
5   & \multicolumn{1}{c|}{5}  & \multicolumn{1}{c|}{1} & \multicolumn{1}{c|}{1} & 1                      & \multicolumn{1}{c|}{0} & 1                      & \multicolumn{1}{c|}{0} & 1                      & 0     \\
6   & \multicolumn{1}{c|}{9}  & \multicolumn{1}{c|}{2} & \multicolumn{1}{c|}{1} & \multicolumn{1}{c|}{1} & \multicolumn{1}{c|}{1} & \multicolumn{1}{c|}{1} & \multicolumn{1}{c|}{1} & \multicolumn{1}{c|}{1} & 1     \\
7   & \multicolumn{1}{c|}{14} & \multicolumn{1}{c|}{2} & \multicolumn{1}{c|}{2} & \multicolumn{1}{c|}{2} & \multicolumn{1}{c|}{2} & \multicolumn{1}{c|}{2} & \multicolumn{1}{c|}{2} & \multicolumn{1}{c|}{1} & 1    
\end{tabular}
\end{table}

We observe that the proposed algorithm does not formally guarantee that the selected points satisfy Theorem~\ref{theor:ChungYao}. To address this issue, one can always perform a test like the one presented in Section~\ref{sec:test1} to assess the extent to which polynomials are included in the local space. If the results are not satisfactory, the number of promoted degrees of freedom can be increased to ensure that all required polynomials are included. In the limit, our method coincides with a standard Lagrange Finite Element Method on the triangulation $\bigcup_{E \in \Th} \TE[E]$. We would like to emphasize that the numerical experiments show that this heuristic procedure consistently achieves high accuracy with the prescribed initial number of degrees of freedom stated in \eqref{eq:totalnumberdofs}.

\subsection{The determination of shape functions}\label{sec:wights}

Given a set of fine nodes, divided into coarse and virtual nodes, we can write new shape functions for the polygon $E$ as a weighted combination of finite element basis functions on $\TE[E]$. More precisely, we write the final shape functions as 
\begin{equation}
    \varphi_{i}(\xx) = \Psi_{i}(\xx) + \sum_{j \in \virtualnodes[E]} \omega_{ij} \Psi_{j}(\xx)\quad \forall i =1,\dots,\Ndof[E], \quad \forall \xx \in E,
    \label{eq:basisfunctions}
\end{equation}
for a given set of weights $\omega_{ij} \in \R$, and for all $i \in \coarsenodes[E]$ and $j \in \virtualnodes[E]$.

We observe that the functions in \eqref{eq:basisfunctions} satisfy Property \ref{prop:krnocker}. Indeed, 
\begin{align*}
    \dof_{\ell}(\varphi_{i}) &= \varphi_{i}(\xx_{\ell}) \\
    &= \Psi_i(\xx_{\ell}) + \sum_{j \in \virtualnodes[E]} \omega_{ij} \Psi_j(\xx_{\ell}) \\
    &= \delta_{i \ell} + \sum_{j = \Ndof[E] + 1}^{N_{\Psi,k}^E} \omega_{ij} \delta_{j\ell} = \delta_{i \ell} \quad\quad \forall i,\ell=1,\dots,\Ndof[E],
\end{align*}
due to Kronecker-Delta property of finite element basis functions with respect to the set of fine nodes. Moreover, Property \ref{prop:continuity} is satisfied by construction, including all nodal coordinates belonging to $\partial E$ into the set of coarse nodes and noticing that $\varphi_{i|\partial E} = \Psi_{i|\partial E}$ for all $i \in \coarsenodes[E]$. Thus, the shape functions $\{\varphi_i\}_{i=1}^{\Ndof[E]}$ satisfy Property \ref{prop:krnocker} and \ref{prop:continuity} independently of the choice of the weights.
Accordingly, these weights can be determined in such a way that the Property \ref{prop:polynomial} is satisfied.

Let us introduce the set of two-dimensional (scaled) monomials of degree up to $k$ on a polygon $E$:
\begin{equation}
    \M{k}{E} = \left\{m_{\alpha}(\xx) = \left(\frac{\xx - \xx_E}{h_E}\right)^{\bm{\alpha}}: \bm{\alpha} = \mathbf{idx}(\alpha)\ \, \forall \alpha = 1,\dots,n_k \right\},
\end{equation}
where $h_E$ denotes the diameter of the polygon $E$ and $\mathbf{idx} : \{1,\dots, n_k\} \to \N^2$ is defined as
\begin{equation}
    \mathbf{idx}(1) = (0,0),\quad \mathbf{idx}(2) = (1, 0),\quad \mathbf{idx}(3) = (0,1),\quad \mathbf{idx}(4) = (2, 0),\quad \dots.
\end{equation}

Since $\{\varphi_i\}_{i=1}^{\Ndof[E]}$ is the set of basis functions of $\VPh[E]{k}$, if a polynomial $p \in \Poly{k}{E}$ belongs to $\VPh[E]{k}$ then it can be expressed as
\begin{equation}
    p(\xx) = \sum_{i=1}^{\Ndof[E]} \dof_i(p) \varphi_i(\xx) \quad \forall \xx \in E.
\end{equation}
Therefore, to ensure that $p\in\VPh[E]{k}$ for any $p\in\Poly{k}{E}$, we can choose the weights $\omega_{ij}$ such that the following equations are satisfied for each $\alpha = 1,\dots,n_k$ and each virtual node $\pp_n$, with $n \in \virtualnodes[E]$:
\begin{align}
    \nonumber m_{\alpha} (\pp_{n}) &= \sum_{i=1}^{\Ndof[E]} \dof_i(m_{\alpha}) \varphi_i(\pp_n)  \\
    \nonumber&= \sum_{i=1}^{\Ndof[E]} \dof_i(m_{\alpha}) [\Psi_i(\pp_n) + \sum_{j \in \virtualnodes[E]} \omega_{ij} \Psi_j(\pp_{n})]\\
    &= \sum_{i=1}^{\Ndof[E]} \dof_i(m_{\alpha}) \omega_{in},\label{eq:constraint_elementwise}
\end{align}
where the last result derives once more from the Kronecker-Delta property of the finite element basis functions with respect to the fine nodes.

Let us introduce the matrices $\mathbf{D} \in \R^{n_k \times \Ndof[E]}$, $\mathbf{V} \in \R^{n_k \times \Nvn[E]}$, and $\mathbf{W} \in \R^{\Ndof[E] \times \Nvn[E]}$, whose entries are defined as 
\begin{gather*}
    \mathbf{D}(\alpha, i) := \dof_{i}(m_{\alpha}) = m_{\alpha}(\xx_i) \quad \forall \alpha=1,\dots,n_k,\quad \forall i=1,\dots,\Ndof[E], \\
    \mathbf{V}(\alpha, j - \Ndof[E]) := m_{\alpha}(\pp_j) \quad 
    \forall \alpha=1,\dots,n_k,\quad \forall j \in \virtualnodes[E], \\
    \mathbf{W}(i,j - \Ndof[E]) := \omega_{ij} \quad \forall i=1,\dots,\Ndof[E],\quad \forall j \in \virtualnodes[E].
\end{gather*}
The constraint \eqref{eq:constraint_elementwise} can be equivalently rewritten as
\begin{equation}
    \mathbf{D} \mathbf{W}(:, n) = \mathbf{V}(:, n),\quad \forall n = 1, \dots, \Nvn[E],
    \label{eq:single_system}
\end{equation}
where $\mathbf{W}(:, n)$ and $\mathbf{V}(:, n)$ denote the $n$-th column of $\mathbf{W}$ and $\mathbf{V}$, respectively.

If $\Ndof[E] = n_k$, then $\mathbf{D}$ is a square non-singular matrix and the system of equations \eqref{eq:single_system} admits a unique solution gathered in the matrix $\mathbf{W}$ for any choice of coarse nodes that satisfies Theorem \ref{theor:ChungYao}.
However, when $\Nv[E] > 3$, the local number of degrees of freedom $\Ndof[E]$ exceeds $n_k$ in order to satisfy the conformity requirements stated in \ref{prop:continuity}. Consequently, the problem \eqref{eq:single_system} is typically underdetermined. To overcome this issue, we determine the weights by solving the following minimization problem:
\begin{align}
    \label{eq:nuovooptsumofsquared}\min_{\mathbf{W} \in \R^{\Ndof[E] \times \Nvn[E]}}
        &\sum_{i \in \coarsenodes[E]} \sum_{j \in \virtualnodes[E]} \omega_{ij}^2\\
    \text{such that}\quad 
   \label{eq:nuovooptvincolo} &\mathbf{D} \mathbf{W}(:, n) = \mathbf{V}(:, n),\quad \forall n = 1, \dots, \Nvn[E].
\end{align}
We observe that the choice of minimizing the sum of the squared weights means that we look for a uniform contribution of all the virtual functions to determine the coarse functions \cite{Bunge2020}.

To solve the minimization problem \eqref{eq:nuovooptsumofsquared}-\eqref{eq:nuovooptvincolo}, we can solve the corresponding Karush-Kuhn-Tucker (KKT) system:
\begin{equation}
    \begin{bmatrix}
        \mathbf{Q}& \mathbf{A}^T \\
        \mathbf{A} & \mathbf{O} \\
    \end{bmatrix} \begin{bmatrix} \bm{\omega} \\ \bm{\lambda}\end{bmatrix} = \begin{bmatrix}
        \bm{0}\\
        \bm{b}
    \end{bmatrix},
    \label{eq:kkt}
\end{equation}
where 
\begin{equation*}
    \begin{gathered}
        \mathbf{Q} = 2\, \mathbf{I} \in \R^{\Ndof[E]\Nvn[E] \times \Ndof[E]\Nvn[E]},\\
        \mathbf{A}  \in \R^{n_k \Nvn[E] \times \Ndof[E] \Nvn[E]}: \mathbf{A} = \operatorname{diag}(\mathbf{D}^T, \dots, \mathbf{D}^T), \\
        \mathbf{b}  \in \R^{n_k \Nvn[E]}:\quad \mathbf{b} = [\mathbf{V}(:, 1),\dots,\mathbf{V}(:, \Nvn[E])],\\
        \bm{\omega}  \in \R^{\Ndof[E] \Nvn[E] }:\quad \bm{\omega} = [\mathbf{W}(:, 1),\dots,\mathbf{W}(:, \Nvn[E])],
    \end{gathered}
\end{equation*}
$\bm{\lambda} \in \R^{n_k \Nvn[E]}$ denotes the vector of Lagrangian multipliers associated with the constraint \eqref{eq:nuovooptvincolo}, $\mathbf{I}\in\R^{\Ndof[E]\Nvn[E] \times \Ndof[E]\Nvn[E]}$ is the identity matrix, $\mathbf{O} \in \R^{n_k \Nvn[E] \times n_k \Nvn[E]}$ is the zero matrix, and $\bm{0} \in \R^{\Ndof[E]\Nvn[E]}$ is the zero vector.

\begin{proposition}\label{prop:optproblemsolution}
    Given $\Ndof[E] \geq n_k$ and a choice of coarse nodes satisfying Theorem \ref{theor:ChungYao}, Problem \eqref{eq:kkt} admits a unique solution and $\Poly{k}{E} \subset \VPh[E]{k}$. 
\end{proposition}
\begin{proof}
    Observe that the optimization Problem \eqref{eq:kkt} is equivalent to 
\begin{equation*}
    \begin{cases}
        2 \bm{\omega} + \mathbf{A}^T \bm{\lambda} = \bm{0}\\
        \mathbf{A} \bm{\omega} = \bm{b}
    \end{cases} \Leftrightarrow \begin{cases}
        \bm{\omega} = -\frac{1}{2} \mathbf{A}^T \bm{\lambda} \\
        \mathbf{A} \bm{\omega} = \bm{b}
    \end{cases}\Leftrightarrow\begin{cases}
        \bm{\omega} = -\frac{1}{2} \mathbf{A}^T \bm{\lambda} \\
        \mathbf{A} \mathbf{A}^T \bm{\lambda} = -2 \bm{b}
    \end{cases}
\end{equation*}
The matrix $\mathbf{A} \mathbf{A}^T$ is a non-singular diagonal block matrix, since its diagonal blocks are given by the squared non-singular matrix $\mathbf{D}^T \mathbf{D}$. Indeed, thanks to Theorem \ref{theor:ChungYao}, matrix $\mathbf{D}$ is a full-rank matrix by rows, and, as a consequence, $\mathbf{D}^T \mathbf{D}$ is non-singular. Thus, the solution of Problem \eqref{eq:kkt} exists and is unique.
\end{proof}
By computing the Cholesky factorization \cite{higham2009cholesky} of the small matrix $\mathbf{D}^T \mathbf{D}$, it is possible to efficiently solve system \eqref{eq:kkt}, as detailed in Algorithm \ref{alg:kkt}.

\SetKwComment{Comment}{/* }{ */}

\begin{algorithm}
\caption{An algorithm to efficiently solve system \eqref{eq:kkt}. }\label{alg:kkt}
\KwData{The matrices $\mathbf{D}$ and $\mathbf{V}$.}
\KwResult{The shape function weights $\mathbf{W}$.}
$\mathbf{L} \gets \operatorname{chol}(\mathbf{D}^T\mathbf{D})$ \Comment*[r]{Compute Cholesky factorization of $\mathbf{D}^T\mathbf{D} = \mathbf{L} \mathbf{L}^T$ once.}
\For{$j = 1,\dots, \Nvn[E]$}{ 
    $\mathbf{L} \bm{y} = \mathbf{V}(:, j) \rightarrow \bm{y}$\;
    $\mathbf{L}^T \bm{x} = \bm{y} \rightarrow \bm{x}$\Comment*[r]{Solve two triangular systems of dimension $n_k$ at each iteration.}
    $\mathbf{W}(:, j) \gets \mathbf{D}^T \bm{x}$\;
}
\end{algorithm}

\begin{remark}

For $k=2$, we recall that a similar approach was already proposed in \cite{Bunge2022} to construct second-order shape functions. In that work, an optimization problem with the same cost functional but different constraints was solved to determine the shape function weights. Specifically, for $k=2$, \cite{Bunge2022} considers the following minimization problem:
\begin{align}
    \label{eq:opt:minsumsquared}\min_{\mathbf{W} \in \R^{\Ndof[E] \times \Nvn[E]}}\quad
    &\sum_{i \in \coarsenodes[E]} \sum_{j \in \virtualnodes[E]} \omega_{ij}^2\\
    \label{eq:opt:weightsunity}\text{such that}\quad &\sum_{i \in \coarsenodes[E]} \omega_{ij} = 1\quad \forall j \in \virtualnodes[E]\\
    \label{eq:opt:weightsaffine}&\pp_{j} = \sum_{i \in \coarsenodes[E]} \omega_{ij} \xx_i \quad \forall j \in \virtualnodes[E].
\end{align}
These constraints ensure that the space of linear polynomials is contained in the final approximation space.
For this reason, this construction is not sufficient to reproduce the desired 
convergence rate of $\Poly{2}{}$-elements \cite{Bunge2022}, since it only yields a linear convergence rate.

To overcome this limitation, \cite{Bunge2022} proposed replacing the cost functional \eqref{eq:opt:minsumsquared} with the following variational energy functional:
\begin{equation}
    \sum_{i \in \coarsenodes[E]} \left[\sum_{e \in \Eh[E,\virtualnodes]} \int_e \norm[2]{\nabla^+ \varphi_i^E - \nabla^- \varphi_i^E}^2 + \epsilon \int_E \norm[2]{\nabla \varphi_i^E}^2 \right], 
    \label{eq:opt:varenergy}
\end{equation}
where $\Eh[E,\virtualnodes]$ denotes the set of virtual edges related to the element $E$, and the operators $ \nabla^+_\sigma$ and $ \nabla^+_\sigma$ correspond to the gradient with respect to the triangle on the right ($+$) and on the left ($-$) of the edge, respectively. 
The first term of the cost function \eqref{eq:opt:varenergy} penalizes flux jumps, whereas the second one acts as a regularization term, which is needed since the flux matrix may be singular for non-simplicial cells.

By adopting this formulation, \cite{Bunge2022} achieved the expected convergence rate for 
$k=2$.
It is worth noting that the weights obtained through \eqref{eq:opt:varenergy} generally produce smoother shape functions than ours, typically leading to slightly smaller error constants. However, extending Problem \eqref{eq:opt:varenergy} to higher polynomial orders is nontrivial, as it would require penalizing higher-order derivative jumps and properly balancing them within a single cost functional. Moreover, such an extension entails significantly higher computational cost.
\end{remark}

\section{The Zipped Discrete Problem}\label{sec:discreteproblem}

In the following, we add a superscript $E$ every time we want to highlight the dependency of the element $E$. Moreover, we will use the same symbol $C$ to denote a positive constant that does not depend on the meshsize $h$, with different meanings in different contexts.

Given the shape functions described in the previous section, we are now able to properly describe the final local zipped finite element space introduced in Equation~\eqref{eq:ele_space} as
\begin{align}
    \VPh[E]{k} := \{ v \in \con{0}{\overline{E}} \cap \sob{1}{E}: i) &\ v \in \mathbb{P}_k(T), \quad \forall  T \in \TE[E], \\
    \label{eq:ZemblConstraint} ii) &\  v(\pp_j^E)= \sum_{i \in \coarsenodes[E]} \omega_{ij}^E\dof^{E}_i(v), \quad \forall j \in \virtualnodes[E] \}, 
\end{align}
where the triangulation $\TE[E]$ is defined as in \eqref{eq:triangles}. We remark that $\Poly{k}{E} \subset \VPh[E]{k}$ due to Proposition \ref{prop:optproblemsolution}.

On $\VPh[E]{k}$, we choose the following set of degrees of freedom:
\begin{enumerate}[label={\textbf{D.\arabic*}}]
    \item \label{dof:vertices}values of $v$ at the vertices of $E$;
    \item \label{dof:edge} if $k\geq2$, the values in $k-1$ evenly spaced points internal of the edges of $e \in \Eh[E]$; 
    \item \label{dof:internal} if $k\geq3$, a set of internal nodes that allows to satisfy Theorem \ref{theor:ChungYao}.
\end{enumerate}

\bigskip

\begin{proposition}
    The degrees of freedom \ref{dof:vertices}-\ref{dof:internal} are unisolvent for $\VPh[E]{k}$.
\end{proposition}
\begin{proof}
We observe that the number of degrees of freedom is equal to the dimension of the space. 
To prove the unisolvence, we have to show that
\begin{equation*}
    \dof_i^E(v) = 0 \quad \forall i\in \coarsenodes[E] \quad\Rightarrow \quad v = 0,\quad \forall v \in \VPh[E]{k}.
\end{equation*}

Thus, from the DOFs definition of Equation~\eqref{eq:dofs}, let us consider a generic function $v\in \VPh[E]{k}$ such that
\begin{equation}
     v(\bm{\xi}_i) = 0 \quad \forall i \in \coarsenodes[E].
     \label{eq:zerodofs_coarse}
\end{equation}

Moreover, from Equations \eqref{eq:ZemblConstraint} and \eqref{eq:zerodofs_coarse}, it holds
\begin{equation}
    v(\bm{\xi}_i) = 0 \quad \forall i \in \virtualnodes[E].
    \label{eq:zerdofs_virtual}
\end{equation}

Defined $\mathbb{V}(E; \TE[E])$ as in \eqref{eq:fintielementlocalspace}, we observe that the evaluations of a function $v \in \VPh[E]{k} \subset \mathbb{V}(E; \TE[E])$ at the fine nodes, i.e. $v(\bm{\xi}_i),\ \forall i \in \finenodes[E]$ constitute a set of unisolvent degrees of freedom for the space $\mathbb{V}(E; \TE[E])$.
Since $\finenodes[E] = \coarsenodes[E] \cup \virtualnodes[E]$, Equations \eqref{eq:zerodofs_coarse} and \eqref{eq:zerdofs_virtual} imply $v = 0$.

\end{proof}

Given the global Z-FEM space defined in \eqref{eq:global_space}, the Zipped discretization of Problem \eqref{eq:varproblem} reads as: \textit{Find $u_h \in \VPh{h, k}$ such that} 
\begin{equation}
    \bilin{u_h}{v_h} = \scal[\Omega]{f}{v_h} \quad \forall v_h \in \VPh{h, k}.
    \label{eq:discreteproblem}
\end{equation}

As in the standard FEM formulation, the existence and uniqueness of the solution of Problem \eqref{eq:discreteproblem} are inherited by the continuous framework. Moreover, the continuity and the coercivity of the bilinear form $\mathcal{B}$ on $\VPh{h,k} \subset \VP$ are the key ingredients for C\'ea's Lemma \cite{brennerscott}:
    \begin{equation}
    \label{eq:CeaLemma}
        \norm[\sob{1}{\Omega}]{u-u_h} \leq C \min_{v_h \in \VPh{h,k}} \norm[\sob{1}{\Omega}]{u-v_h},
    \end{equation}
where $C$ is a constant related to continuity and coercivity constants of $\mathcal{B}$ on $\VP$.

\begin{proposition}
\label{prop:stimeapriori}
    Under Mesh Assumption \ref{ass:meshassumption}, let $u \in \sob{k+1}{\Omega}$, $k \geq 1$, be the solution of \eqref{eq:varproblem}, and let $u_h \in \VPh{h, k}$ be the solution of the discrete problem \eqref{eq:discreteproblem}. It holds:
    \begin{equation}
         \norm[\sob{1}{\Omega}]{u-u_h} \leq Ch^{k} |u|_{\sob{k+1}{\Omega}},
         \label{eq:apriori1}
    \end{equation}
    and
    \begin{equation}
         \norm[\leb{2}{\Omega}]{u-u_h} \leq Ch^{k+1} |u|_{\sob{k+1}{\Omega}}.
         \label{eq:apriori2}
    \end{equation}
\end{proposition}

\begin{proof}
Let us consider an element $E \in \Th$. Let us define the interpolation operator $\mathcal{I}_k^E : \sob{1}{E} \to \VPh[E]{k}$ for all $E \in \Th$:

\begin{equation}
    \mathcal{I}_k^E v := \sum_{i \in \coarsenodes[E]} \dof_i(v) \varphi_i. 
    \label{eq:interpE}
\end{equation}

We define the set of coarse nodes that belong to triangle $T \in \TE[E]$ as $\mathcal{C}^E_T \subset \coarsenodes[E]$
 and the virtual nodes that belong to triangle $T \in \TE[E]$ as $\virtualnodes[E]_T \subset \virtualnodes[E]$.
The interpolation operator $ \mathcal{I}_k^T: \sob{1}{T} \to \Poly{k}{T}$ on $T \in \TE[E]$ is defined as the restriction of \eqref{eq:interpE} on $T$:
\begin{align*}
    \mathcal{I}_k^T v &:= \sum_{i \in \coarsenodes[E]_T} \dof_i(v) \varphi_i|_T \\
                      &= \sum_{i \in \coarsenodes[E]_T} \dof_i(v) [\Psi_i|_T +\sum_{j \in \virtualnodes[E]_T} \omega_{ij} \Psi_j|_T]  \\
                      &= \sum_{j \in \coarsenodes[E]_T \cup \virtualnodes[E]_T} \overline{\dof}_j(v) \Psi_j|_T, \quad \quad \forall v \in \sob{1}{T},
\end{align*}
where $ \{ \overline{\dof}_j(v) \}_{j \in\coarsenodes[E]_T \cup \virtualnodes[E]_T} $ are defined as: $\forall v \in \sob{1}{T}$
\begin{equation*}
    \overline{\dof}_i(v) = \left\{
    \begin{aligned}
        &\dof_i(v), && i \in \coarsenodes[E]_T,\\
        &\sum_{\ell \in \coarsenodes[E]_T} \omega_{\ell i} \dof_{\ell}(v), &\quad& i \in \virtualnodes[E]_T.
    \end{aligned}\right.
\end{equation*}

Note that if $v \in \sob{k+1}{T}$, then $\mathcal{I}_k^T v  \in \sob{k+1}{T}$. Moreover, since $\Poly{k}{E} \subseteq \VPh[E]{k}$, then $p=\mathcal{I}_k^E p$ and $p_{|T} = \left(\mathcal{I}_k^Ep \right) _{|T} = \mathcal{I}_k^T p$. Thus, from Theorem 4.4.4 in \cite{brennerscott}, it follows that: 
\begin{equation*}
    |u - \mathcal{I}_k^{T} u|_{\sob{1}{T}} \leq C h_T^{k}  |u|_{\sob{k+1}{T}}.
\end{equation*}

The Mesh Assumptions \ref{ass:meshassumption} imply that the triangles belonging to $\TE[E]$ of an element $E \in \Th$ are regular (see Remark~\ref{rem:shaperegular}).

Now, by Theorem 4.4.20 in \cite{brennerscott}, we derive the following estimate on $E \in \Th$, since $(u - \mathcal{I}_k^E u) \in \sob{1}{E}$,
\begin{equation*}
    \norm[\sob{1}{E}]{u - \mathcal{I}_k^E u} = \left(\sum_{T \in \TE[E]}  \norm[\sob{1}{T}]{u - \mathcal{I}_k^T u}^2 \right)^{1/2} \leq C \left(\max_{T \in \TE[E]}{h_T}\right)^{k} |u|_{\sob{k + 1}{E}} \leq C h_E^k |u|_{\sob{k + 1}{E}}.
\end{equation*}
Thus, 
\begin{equation}
\label{eq:interpglobalestimate}
    \norm[\sob{1}{\Omega}]{u - \mathcal{I}_k u} \leq C h^k |u|_{\sob{k+1}{\Omega}},
\end{equation}
where $(\mathcal{I}_k)_{|E} = \mathcal{I}_k^E$ for all $E \in \Th$.

From C\'ea's Lemma \eqref{eq:CeaLemma}, we have
\begin{align}
\label{eq:proofH1final}
    \norm[\sob{1}{\Omega}]{u-u_h} \leq C \min_{v_h \in \VPh{h,k}} \norm[\sob{1}{\Omega}]{u-v_h} \leq C\norm[\sob{1}{\Omega}]{u - \mathcal{I}_k u} \leq C h^k |u|_{\sob{k + 1}{E}} .
\end{align}
Equation \eqref{eq:apriori1} follows.

The $L^2$-norm estimate follows from a duality argument. Let us define the adjoint operator related to the problem \eqref{eq:modelproblem} as $\mathcal{L}^*w := -\nabla \cdot (\D \nabla w) + \gamma w$. For every $g \in \leb{2}{\Omega}$, there exists a unique $w \in \sob{2}{\Omega} \cap \VP $ such that 
\begin{equation}
\label{eq:adjointprob}
    \mathcal{L}^*w = g
\end{equation} and 
\begin{equation}
\label{eq:adjointproperty}
    \norm[\sob{2}{\Omega}]{w} \leq C \norm[L^2(\Omega)]{g},
\end{equation}
with a constant $C$ independent of $g$. Let $\phi$ be the solution of the adjoint problem \eqref{eq:adjointprob} associated with the right-hand side $g=u-u_h$. It holds,
\begin{align*}
\|u - u_h\|_{L^2(\Omega)}^2
 &=\scal[\Omega]{u-u_h}{g}= \mathcal{B}(u - u_h, \phi) & \\
 &= \mathcal{B}(u - u_h, \phi - \mathcal{I}_k \phi)\quad \quad &\text{(using Galerkin orthogonality: $\mathcal{B}(u - u_h, \mathcal{I}_k \phi) = 0$ )} \\
 &\le C \norm[\sob{1}{\Omega}]{u - u_h} \norm[\sob{1}{\Omega}]{\phi - \mathcal{I}_k \phi} \quad \quad &\text{(using continuity of $\mathcal{B}$)} \\
 &\le C h \norm[\sob{1}{\Omega}]{u - u_h} \|\phi\|_{H^2(\Omega)} \quad \quad &\text{(using \eqref{eq:interpglobalestimate})} \\
 & \le C h \norm[\sob{1}{\Omega}]{u - u_h}  \norm[L^2(\Omega)]{u-u_h} \quad \quad &\text{(using \eqref{eq:adjointproperty})}.
\end{align*}
By simplifying the previous equation and using \eqref{eq:proofH1final}, the thesis follows.

\end{proof}

\section{Numerical Experiments}
\label{sec:numericalexperiments}
In this section, we show various numerical experiments to assess the flexibility and performance of the Zipped finite element method. 
In Section \ref{sec:test1} we prove that the polynomials can be effectively represented as linear combinations of the Z-FEM basis functions, whereas in Section \ref{sec:test2} we solve Problem~\eqref{eq:modelproblem} and compare the method with FEM and VEM on general meshes.

\subsection{Test 1: Polynomials Approximation}\label{sec:test1}
\begin{table}[!ht]
\centering
  \caption{Test 1: Polynomial approximation errors \eqref{eq:polapproxerrors} over different polygons for different orders.} 
  \label{tab:chap6:singularvalues}
  \begin{tabular}
      {ccc}
      \includegraphics[width=1in]{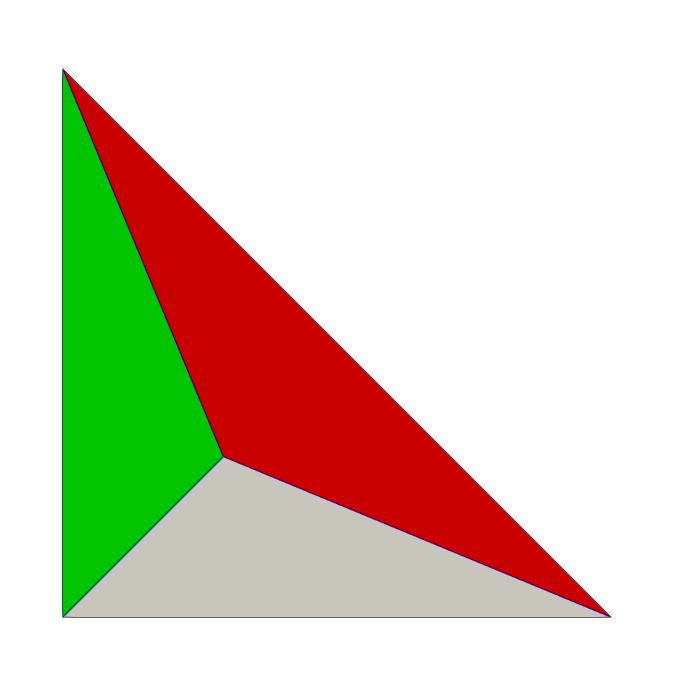} &
      \includegraphics[width=1in]{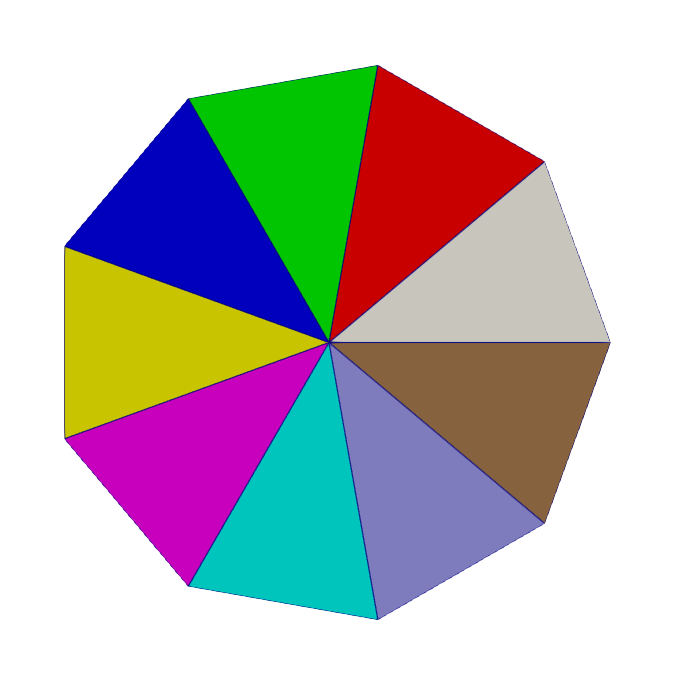} &
      \includegraphics[width=1in]{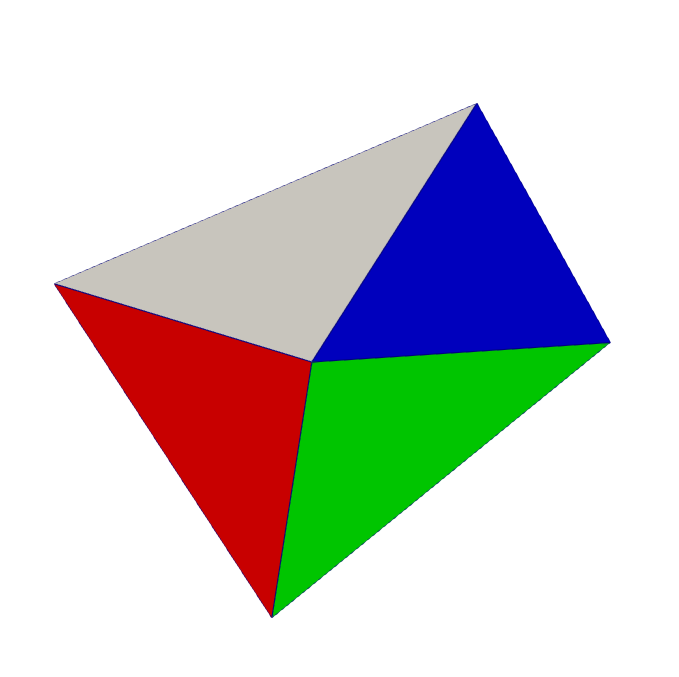}\\
      \multicolumn{1}{c}{Triangle} & \multicolumn{1}{c}{Regular} & \multicolumn{1}{c}{Irregular} \\
      $\Nv[E] = 3$ & $\Nv[E] = 9$ & $\Nv[E] = 4$ \\
\begin{tabular}{ccc}
$k$ & $\mathrm{err}_{I,0}$ & $\mathrm{err}_{I,\nabla}$ \\
1   & 1.8326e-17           & 5.7941e-16                \\
2   & 2.8347e-16           & 2.9336e-15                \\
3   & 1.2973e-15           & 1.8339e-14                \\
4   & 5.8797e-14           & 7.3820e-13                \\
5   & 9.1859e-14           & 3.3798e-12                \\
6   & 5.6793e-13           & 1.4408e-11               
\end{tabular} &
      \begin{tabular}{ccc}
$k$ & $\mathrm{err}_{I,0}$ & $\mathrm{err}_{I,\nabla}$ \\
1   & 2.6949e-16           & 2.8693e-16                \\
2   & 1.3028e-14           & 8.5565e-14                \\
3   & 1.1445e-14           & 8.9299e-14                \\
4   & 1.8168e-14           & 2.6553e-13                \\
5   & 4.7007e-14           & 5.7959e-13                \\
6   & 1.5337e-13           & 2.1020e-12               
\end{tabular} & 
      \begin{tabular}{ccc}
$k$ & $\mathrm{err}_{I,0}$ & $\mathrm{err}_{I,\nabla}$ \\
1   & 6.3167e-17           & 2.3435e-16                \\
2   & 8.2278e-16           & 4.2732e-15                \\
3   & 8.9258e-16           & 6.2626e-15                \\
4   & 9.6685e-15           & 9.3015e-14                \\
5   & 3.0512e-14           & 4.0792e-13                \\
6   & 1.9029e-13           & 4.1390e-12               
\end{tabular} \\
      \includegraphics[width=1in]{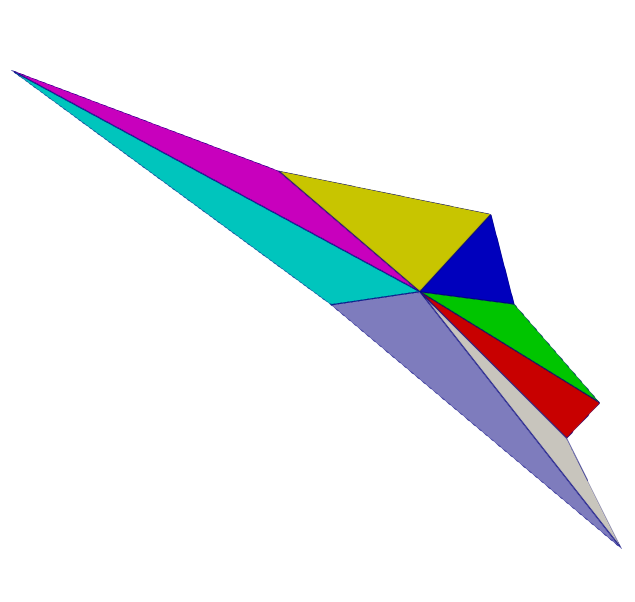} &
      \includegraphics[width=1in]{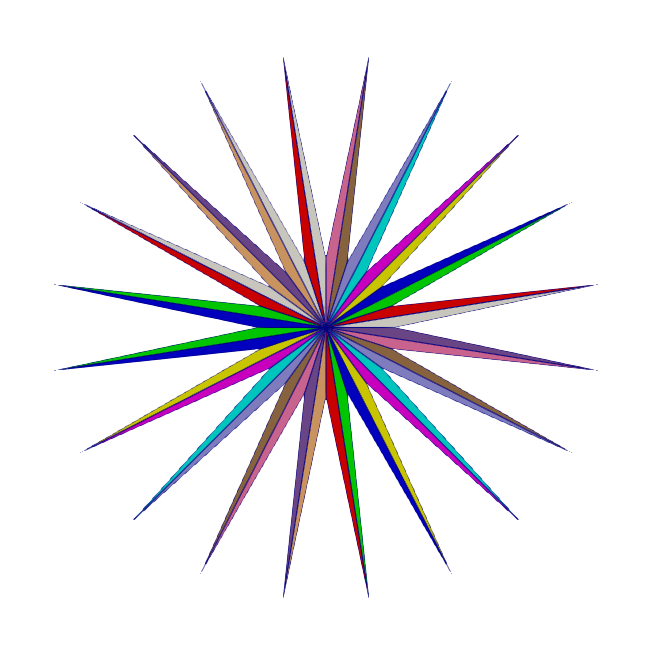} &
      \includegraphics[width=1in]{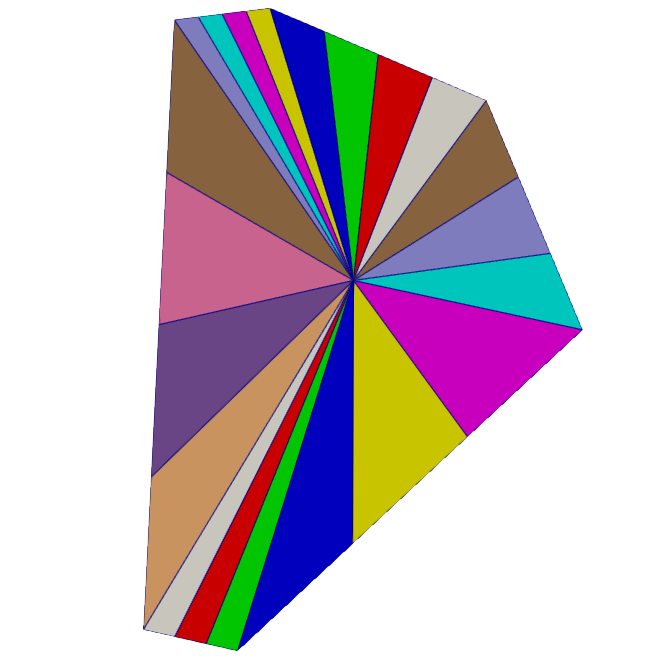}\\
      \multicolumn{1}{c}{Concave} & \multicolumn{1}{c}{Star} &\multicolumn{1}{c}{\begin{tabular}[c]{@{}c@{}}Hexagon\\ with hanging nodes\end{tabular}}  \\
      $\Nv[E]= 8$ & 
      $\Nv[E] = 40$ 
      & $\Nv[E] = 21$  \\
      \begin{tabular}{ccc}
$k$ & $\mathrm{err}_{I,0}$ & $\mathrm{err}_{I,\nabla}$ \\
1   & 1.8493e-17           & 3.0149e-16                \\
2   & 2.1022e-16           & 1.2040e-14                \\
3   & 3.6450e-15           & 2.0565e-13                \\
4   & 1.6660e-13           & 1.9650e-11                \\
5   & 1.3747e-12           & 1.2459e-10                \\
6   & 3.9631e-11           & 1.0054e-08               
\end{tabular} &
      \begin{tabular}{ccc}
$k$ & $\mathrm{err}_{I,0}$ & $\mathrm{err}_{I,\nabla}$ \\
1   & 2.8739e-16           & 6.1237e-15                \\
2   & 2.3066e-16           & 1.9867e-14                \\
3   & 2.1459e-16           & 3.6148e-14                \\
4   & 5.6569e-16           & 1.5407e-13                \\
5   & 3.4445e-15           & 5.2178e-13                \\
6   & 7.2247e-15           & 2.1387e-12               
\end{tabular} &
      \begin{tabular}{ccc}
$k$ & $\mathrm{err}_{I,0}$ & $\mathrm{err}_{I,\nabla}$ \\
1   & 7.1691e-17           & 7.1416e-16                \\
2   & 1.6957e-15           & 1.8354e-14                \\
3   & 2.6098e-15           & 1.4942e-13                \\
4   & 3.2337e-14           & 8.2584e-13                \\
5   & 3.3556e-14           & 4.4693e-12                \\
6   & 6.3862e-13           & 4.8157e-11               
\end{tabular} \\
  \end{tabular}
\end{table}

The first test we propose aims to show how accurately the local basis functions are able to reproduce polynomials. For this purpose, we consider a set of different polygons and, on each polygon, we compute the following polynomial approximation errors
\begin{equation}
    \mathrm{err}_{I,0} = \max_{\alpha =1, \dots, n_k} \norm[\leb{2}{E}]{m_{\alpha} - \mathcal{I}_k^E m_{\alpha}},\quad \mathrm{err}_{I,\nabla} = \max_{\alpha =1, \dots, n_k} \norm[\leb{2}{E}]{\nabla m_{\alpha} - \nabla \mathcal{I}_k^E m_{\alpha}},
    \label{eq:polapproxerrors}
\end{equation}
where the interpolator $\mathcal{I}_k^E: \sob{1}{E} \to \VPh[E]{k}$ is defined in \eqref{eq:interpE}.

Examples from the tested cases are shown in Table~\ref{tab:chap6:singularvalues}. For each one, we draw in the top part the considered polygon and its sub-triangulation, represented by assigning a different color to each sub-triangle. The tested examples include polygons sampled from diverse classes: the reference triangle, a regular nonagon, an irregular quadrilateral, an irregular concave polygon, a star, and, finally, a general polygon with hanging nodes.

For polynomial orders $k = 1, \dots, 6$, we compute and report in Table \ref{tab:chap6:singularvalues} both $\mathrm{err}_{I,0}$ and $\mathrm{err}_{I,\nabla}$, as defined in \eqref{eq:polapproxerrors}. In all tested cases, both errors are found to be close to machine precision, demonstrating that the proposed shape functions can approximate polynomials with very high accuracy, despite the use of the heuristic choice of the internal degrees of freedom described in Section~\ref{sec:shapefunctions}.

\subsection{Test 2: Convergence rates}\label{sec:test2}
\begin{figure}[!ht]
    \centering
    \begin{subfigure}{0.24\textwidth}
        \includegraphics[width=1\linewidth]{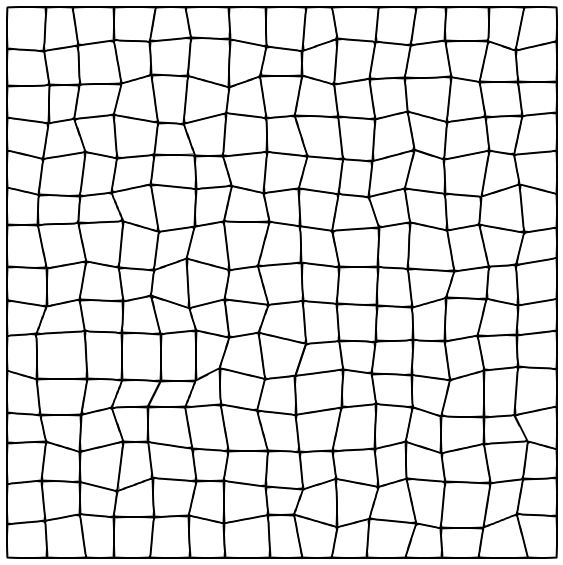}
        \caption{Random Distorted.}
        \label{fig:rdquad_mesh}
    \end{subfigure}
    \begin{subfigure}{0.24\textwidth}
        \includegraphics[width=1\linewidth]{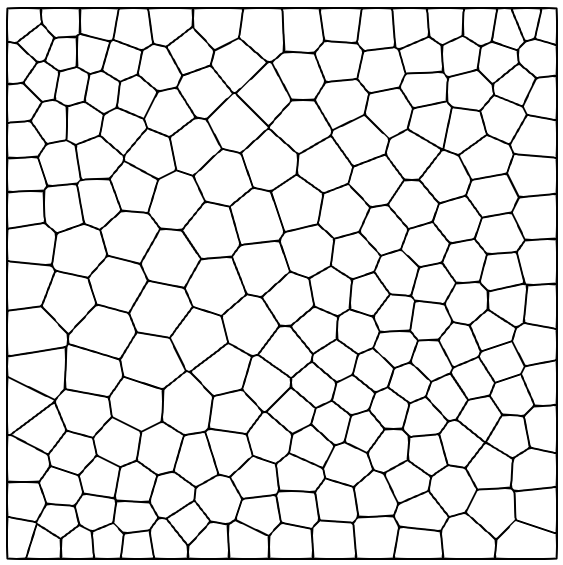}
        \caption{Voronoi.}
        \label{fig:voro_mesh}
    \end{subfigure}
    \begin{subfigure}{0.24\textwidth}
        \includegraphics[width=1\linewidth]{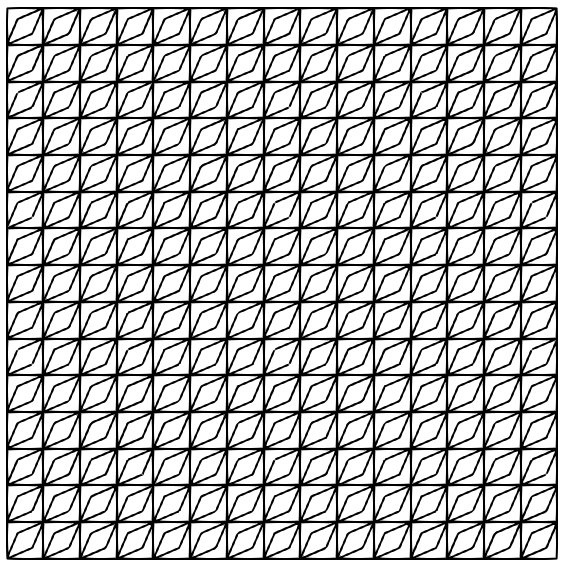}
        \caption{Structured Concave.}
        \label{fig:struconc_mesh}
    \end{subfigure}
    \begin{subfigure}{0.24\textwidth}
        \includegraphics[width=1\linewidth]{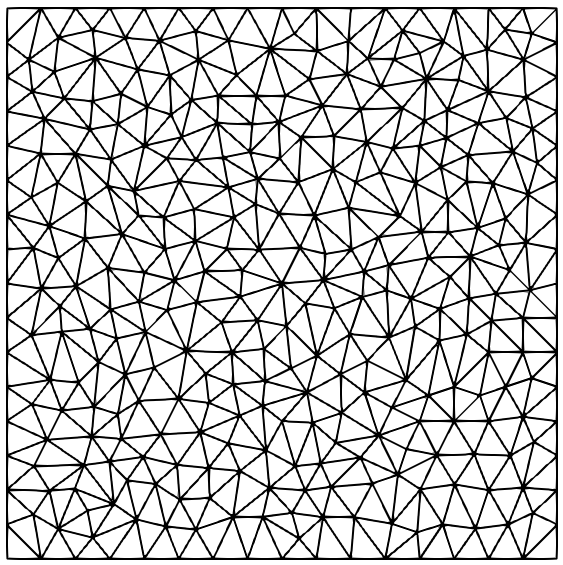}
        \caption{FE triangular mesh.}
        \label{fig:triangle_mesh}
    \end{subfigure}
    \caption{Test 2: Last refinement for each family of meshes employed. }
    \label{fig:test2:mesh}
\end{figure}
%
\begin{table}[!ht]
\caption{Test 2: Empirical order of convergence (EOC) of errors \eqref{eq:errors} computed with respect to $h$, for each value of $k=1,\dots,6$. Each row corresponds to a family of meshes. First: Random Distorted. Second: Voronoi. Third: Structured Concave. The last row refers to Finite Element results computed over a triangular Delaunay family of meshes.}
\label{tab:test2:EOC}
\centering
\begin{tabular}{@{}lccccccc@{}}
\toprule
                     &                         & $k=1$ & $k=2$ & $k=3$ & $k=4$ & $k=5$ & $k=6$ \\ \midrule
\multirow{2}{*}{VEM} & $\mathrm{EOC}_0$        & 1.87  & 3.31  & 4.16  & 5.16  & 6.19  & 7.18  \\
                     & $\mathrm{EOC}_{\nabla}$ & 1.12  & 2.03  & 3.05  & 4.12  & 5.13  & 6.17  \\
\multirow{2}{*}{Z-FEM} & $\mathrm{EOC}_0$        & 1.95  & 3.10  & 4.26  & 5.15  & 6.07  & 7.15  \\
                     & $\mathrm{EOC}_{\nabla}$ & 0.99  & 2.04  & 3.09  & 4.12  & 5.14  & 6.13  \\ \midrule
\multirow{2}{*}{VEM} & $\mathrm{EOC}_0$        & 1.85  & 3.35  & 4.35  & 5.23  & 6.34  & 7.25  \\
                     & $\mathrm{EOC}_{\nabla}$ & 1.10  & 2.10  & 3.20  & 4.24  & 5.38  & 6.34  \\
\multirow{2}{*}{Z-FEM} & $\mathrm{EOC}_0$        & 1.97  & 3.36  & 4.45  & 5.35  & 6.53  & 7.46  \\
                     & $\mathrm{EOC}_{\nabla}$ & 0.98  & 2.20  & 3.30  & 4.31  & 5.46  & 6.41  \\ \midrule
\multirow{2}{*}{VEM} & $\mathrm{EOC}_0$        & 1.67  & 3.08  & 4.01  & 4.91  & 5.97  & 6.81  \\
                     & $\mathrm{EOC}_{\nabla}$ & 1.01  & 1.91  & 2.92  & 3.91  & 4.95  & 5.90  \\
\multirow{2}{*}{Z-FEM} & $\mathrm{EOC}_0$        & 1.77  & 3.00  & 4.07  & 4.89  & 5.97  & 6.90  \\
                     & $\mathrm{EOC}_{\nabla}$ & 0.93  & 1.87  & 2.93  & 3.90  & 4.97  & 5.92  \\ \midrule
\multirow{2}{*}{FEM} & $\mathrm{EOC}_0$        & 1.91  & 3.44  & 4.18  & 5.66  & 6.26  & 7.96  \\
                     & $\mathrm{EOC}_{\nabla}$ & 0.91  & 2.34  & 3.05  & 4.65  & 5.19  & 6.89  \\ \bottomrule
\end{tabular}
\end{table}

Let us consider the Problem~\eqref{eq:modelproblem}, where the diffusion and reaction coefficients are defined as
\begin{equation*}
    \D(\xx) = \begin{bmatrix}
        1 + x_2^2 & -x_1x_2\\
        -x_1x_2 & 1 + x_1^2
    \end{bmatrix},\quad \gamma(\xx) = x_1 x_2,
\end{equation*}
whereas the forcing term and the Dirichlet boundary conditions are such that the exact solution is
\begin{equation*}
    u(\xx) = \sin(2 \pi x_1) \sin(2 \pi x_2).
\end{equation*}

In this test case, we assess the performance of the proposed method by solving this problem on different families of meshes. For each family of meshes, we compute the following standard errors:
\begin{equation}
    \mathrm{err}_0 = \norm[\leb{2}{\Omega}]{u - u_h},\qquad \mathrm{err}_{\nabla} = \norm[\leb{2}{\Omega}]{\nabla u - \nabla u_h},
    \label{eq:errors}
\end{equation}
as the meshsize $h$ decreases (thus, as $\Ndof = \dim \VPh{h,k}$ increases) and for $k=1,\dots, 6$.

Specifically, we consider the following three mesh families:
\begin{itemize}
\item A family of randomly distorted quadrilateral meshes, obtained by perturbing the vertices of standard Cartesian grids;
\item A family of Voronoi meshes generated using \cite{Gedim, Voro};
\item The ``StructuredConcave'' family of meshes available in the PolyDiM repository \cite{Polydim}.
\end{itemize}
The finest meshes of each family are depicted in Figures~\ref{fig:rdquad_mesh}-\ref{fig:struconc_mesh}.

\begin{figure}[!ht]
    \centering
    \begin{subfigure}{0.27\textwidth}
        \includegraphics[width=1\linewidth]{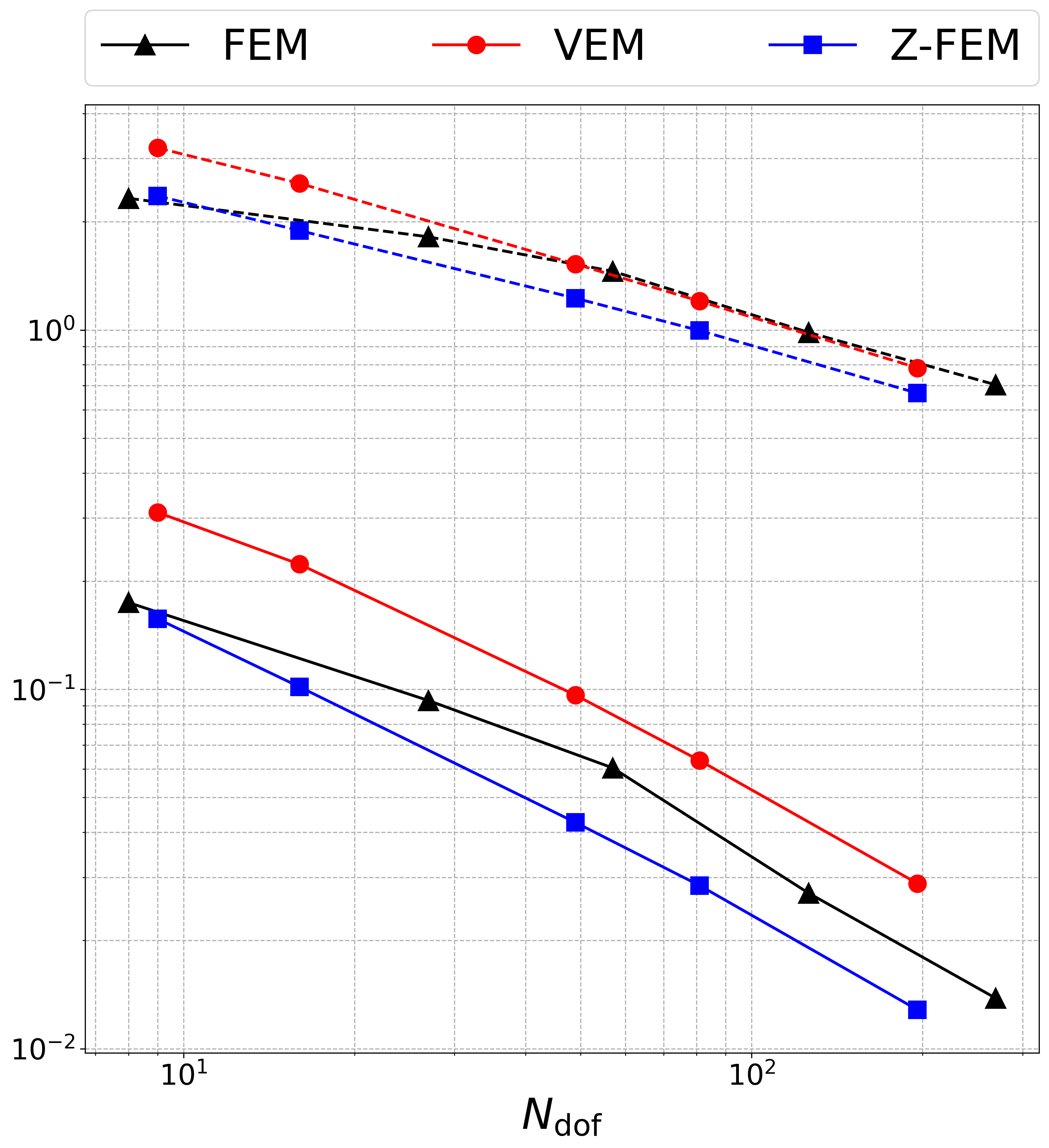}
    \end{subfigure}
    \begin{subfigure}{0.27\textwidth}
        \includegraphics[width=1\linewidth]{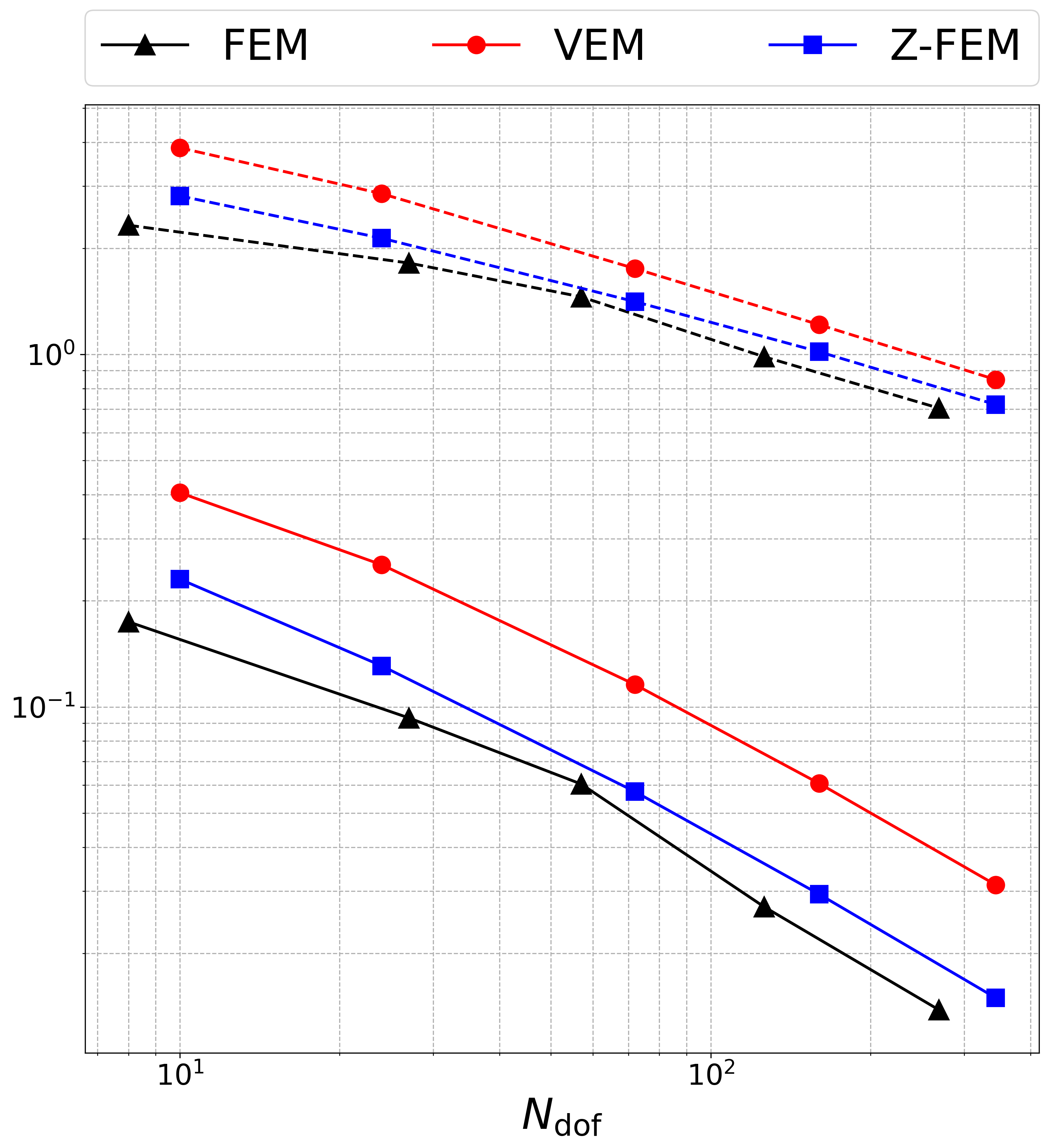}
    \end{subfigure}
    \begin{subfigure}{0.27\textwidth}
        \includegraphics[width=1\linewidth]{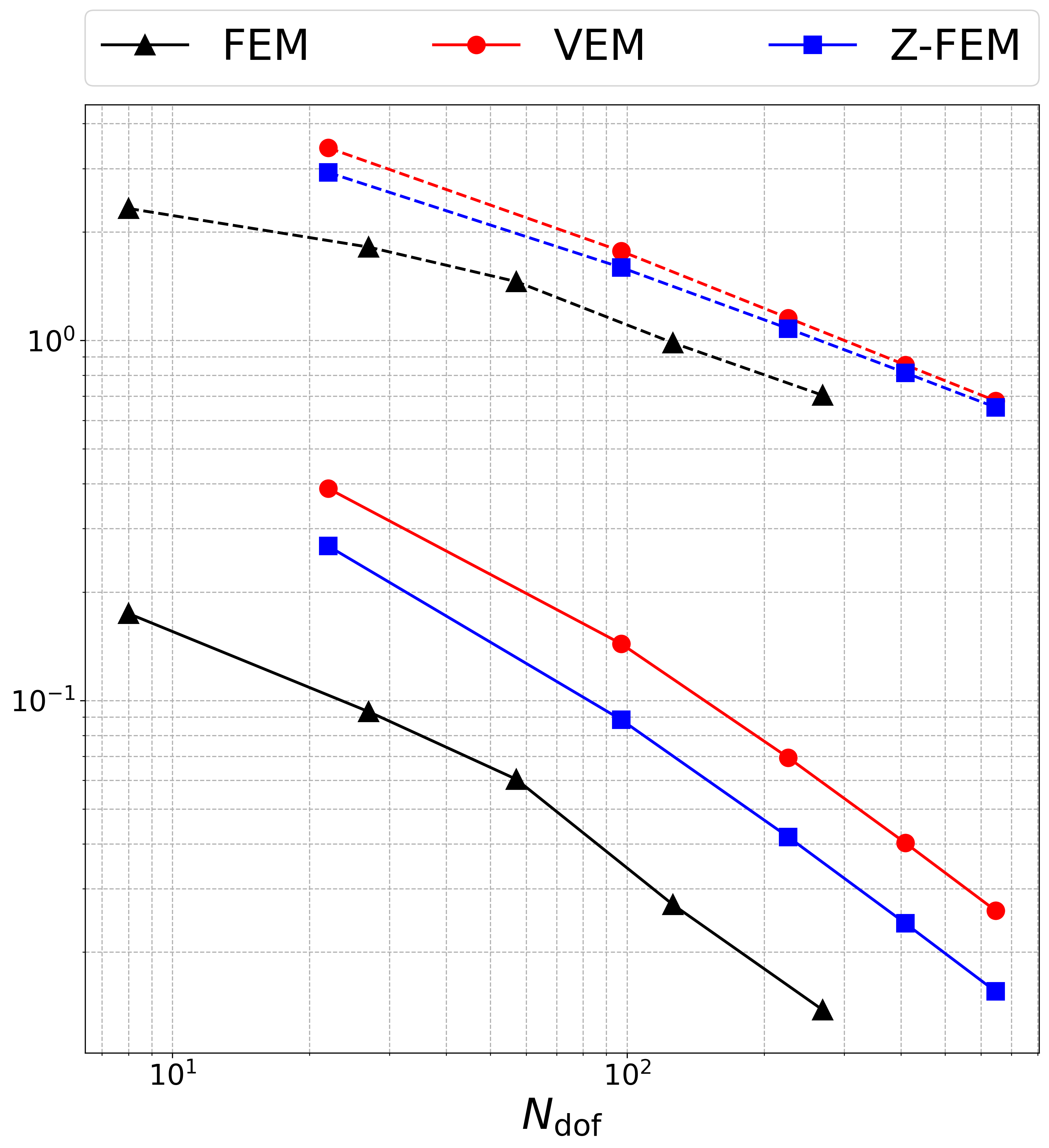}
    \end{subfigure}
    \begin{subfigure}{0.27\textwidth}
        \includegraphics[width=1\linewidth]{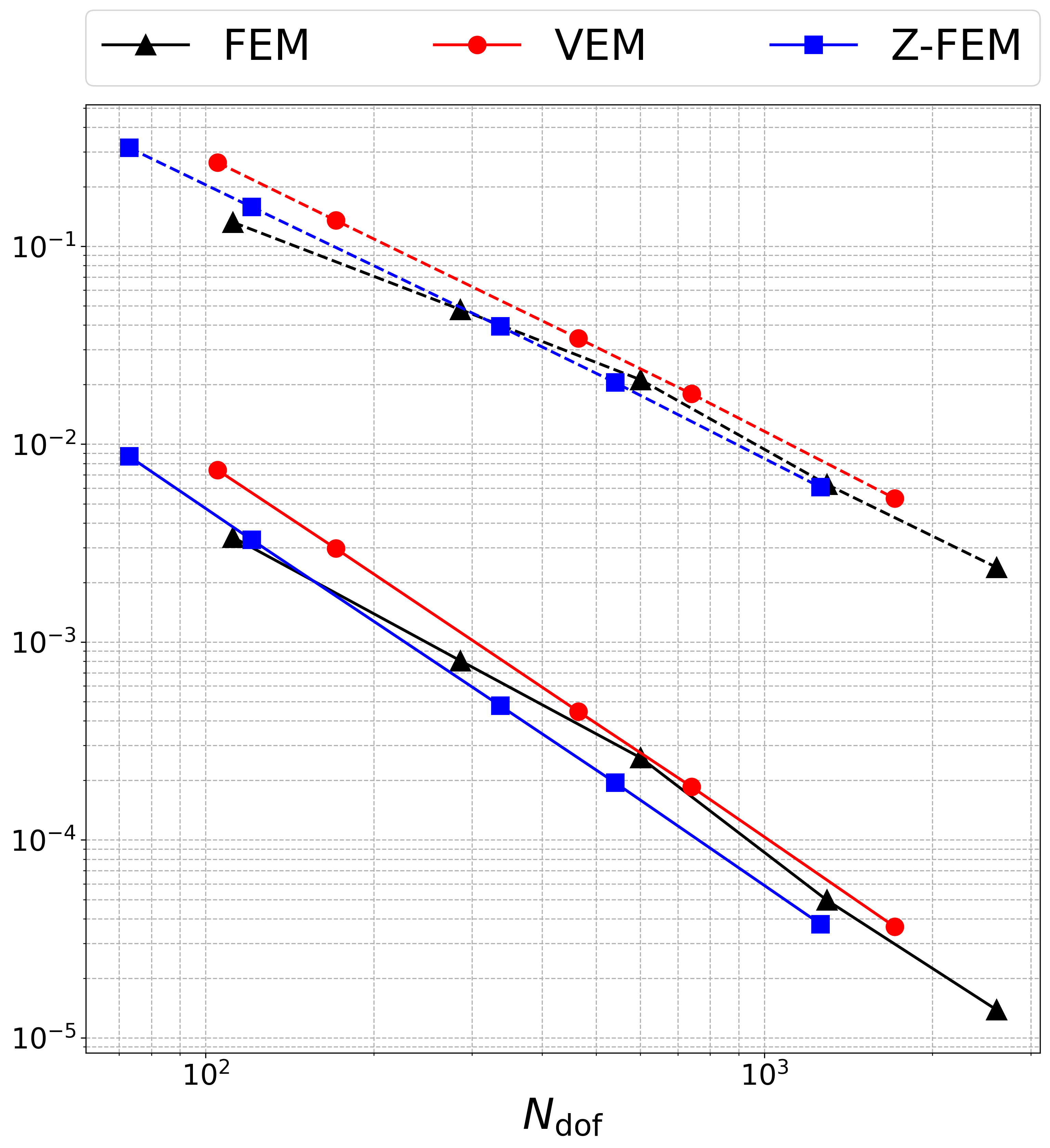}
    \end{subfigure}
    \begin{subfigure}{0.27\textwidth}
        \includegraphics[width=1\linewidth]{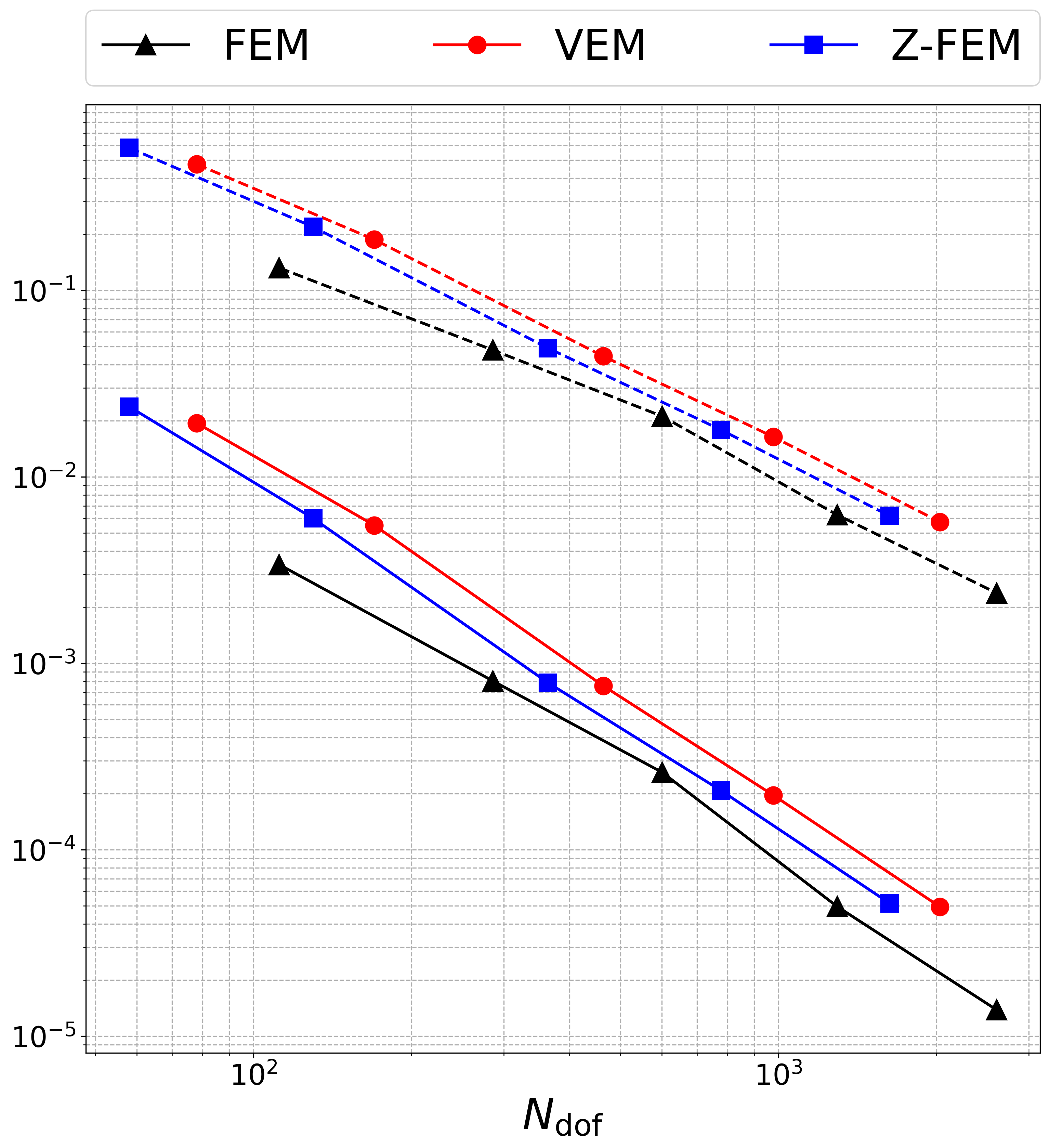}
    \end{subfigure}
    \begin{subfigure}{0.27\textwidth}
        \includegraphics[width=1\linewidth]{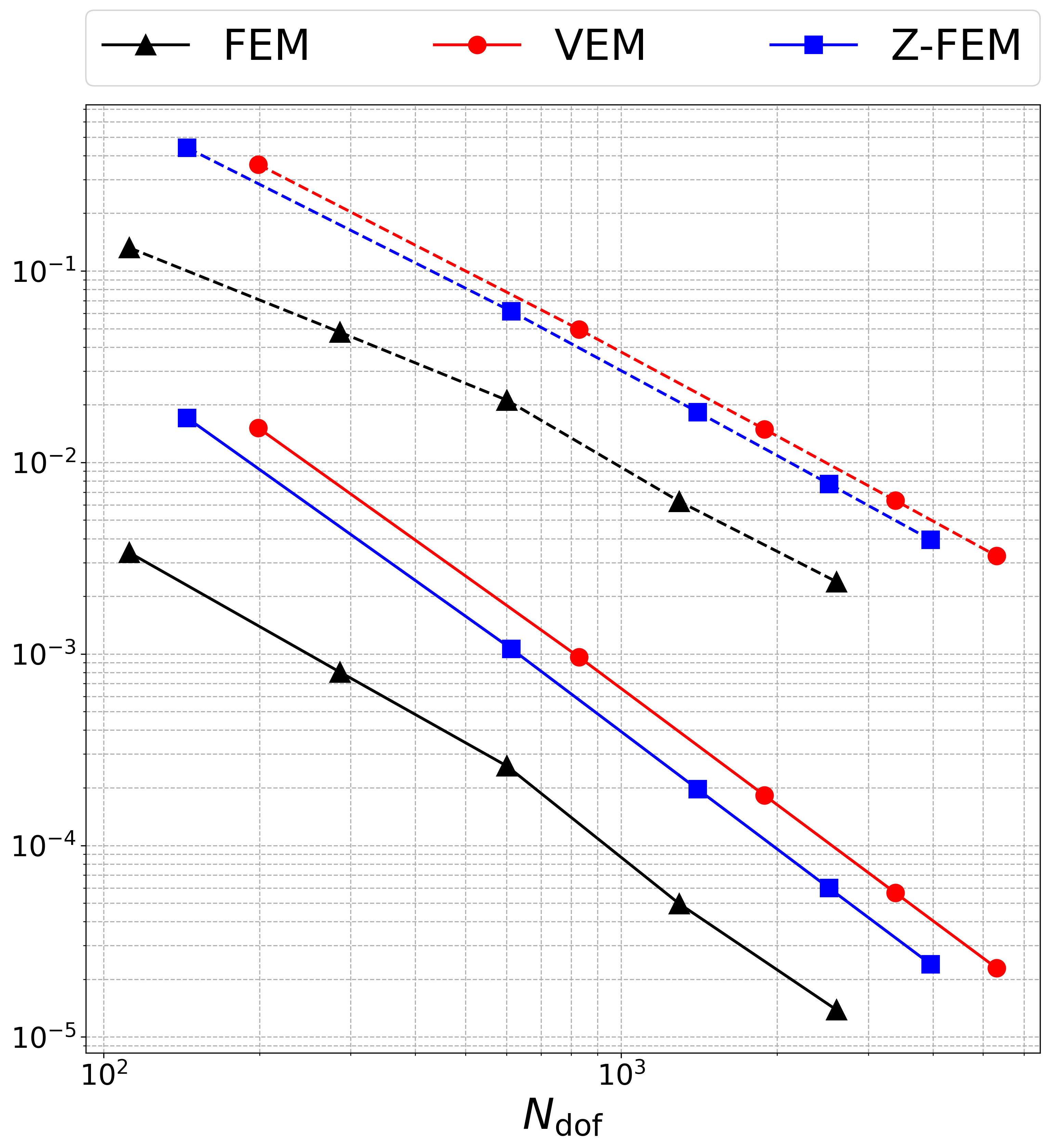}
    \end{subfigure}
    \begin{subfigure}{0.27\textwidth}
        \includegraphics[width=1\linewidth]{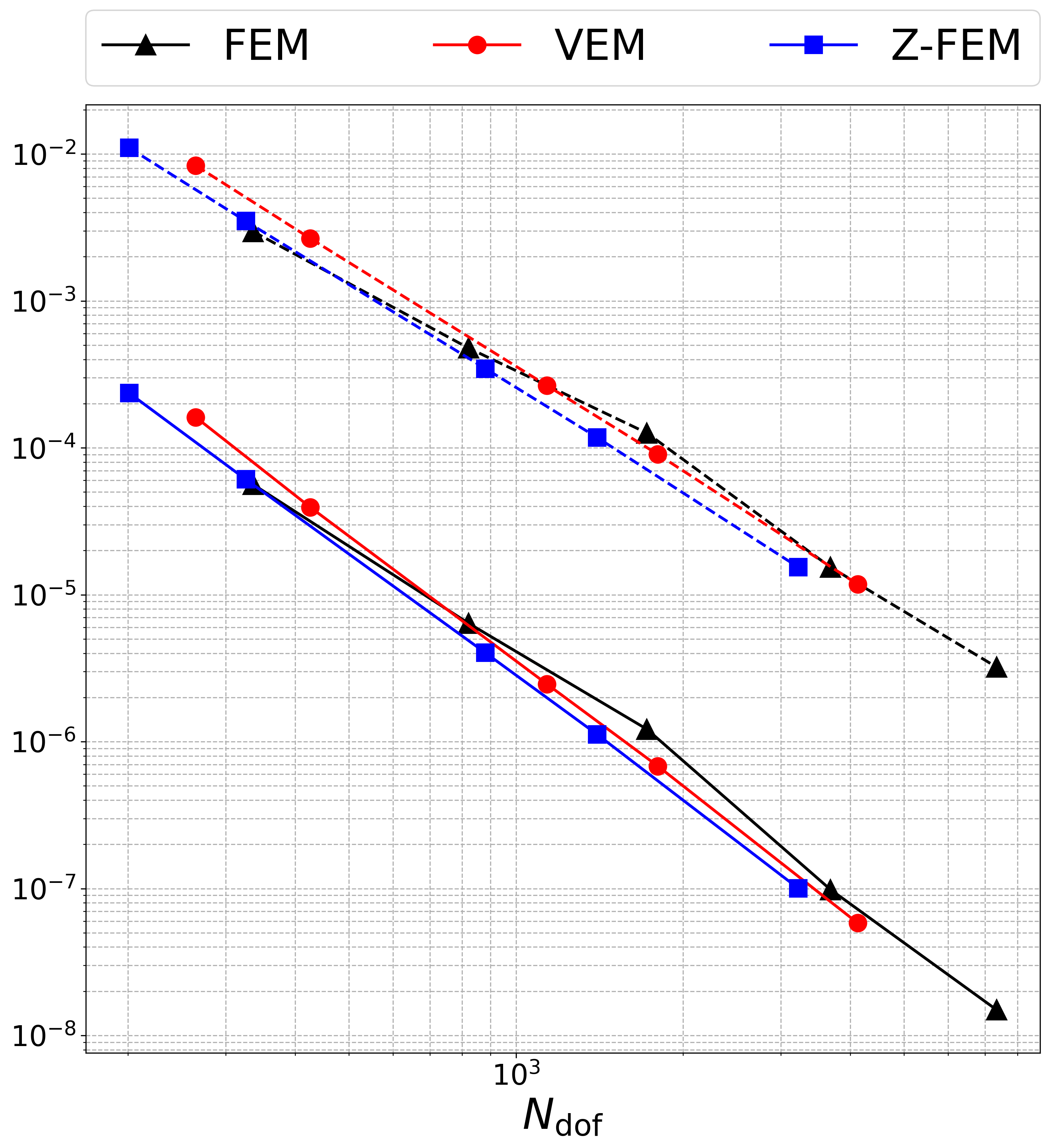}
    \end{subfigure}
    \begin{subfigure}{0.27\textwidth}
        \includegraphics[width=1\linewidth]{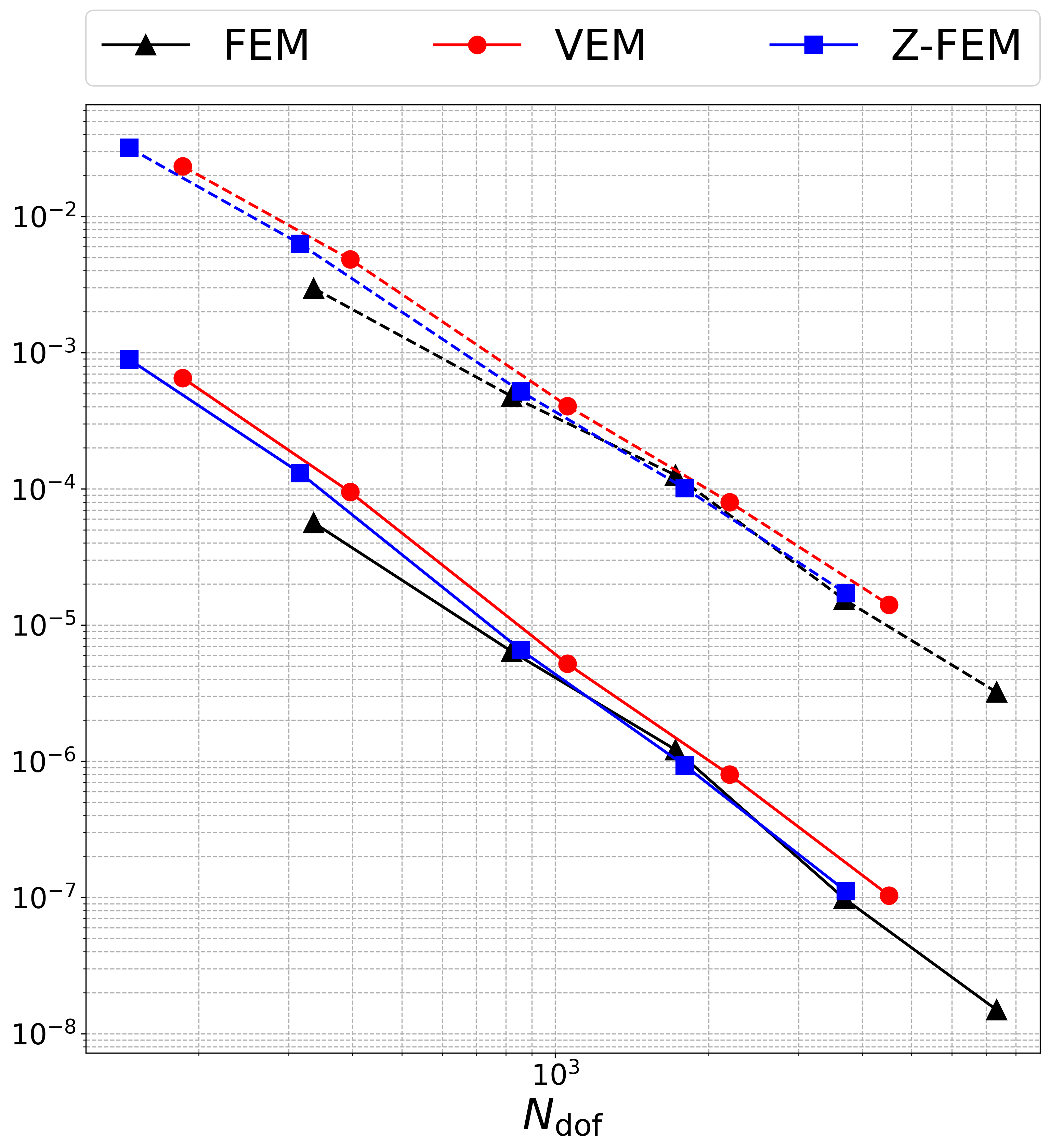}
    \end{subfigure}
    \begin{subfigure}{0.27\textwidth}
        \includegraphics[width=1\linewidth]{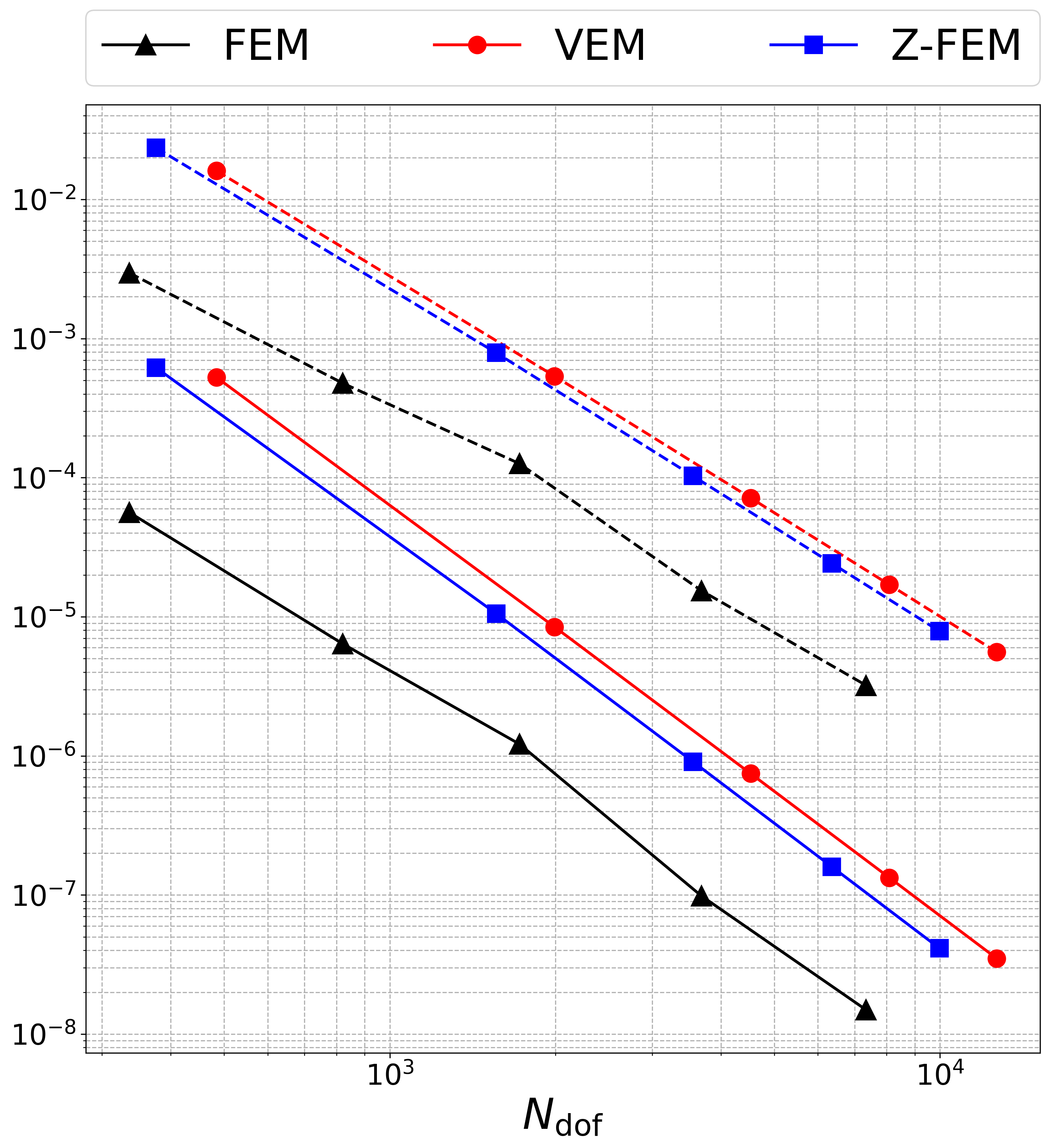}
    \end{subfigure}
    \caption{Test 2: Behavior of errors \eqref{eq:errors} as $\Ndof$ increases. Each column corresponds to a different family of meshes. Left: Random Distorted. Center: Voronoi. Right: Structured Concave. Each row refers to a different value of $k = 1, 3, 5$. Solid lines: $\leb{2}{}$-errors. Dashed lines: $\sob{1}{}$-errors.}
    \label{fig:test2:errors_a}
\end{figure}

In Figure \ref{fig:test2:errors_a}, we report the errors defined in \eqref{eq:errors} as the number of degrees of freedom $\Ndof$ increases, for polynomial orders $k = 1, 3, 5$, and for the three different families of meshes. Moreover, in Table \ref{tab:test2:EOC}, we report the Empirical Order of Convergence for both the $\leb{2}{}$-error (EOC$_0$) and the $\sob{1}{}$-seminorm error (EOC$_\nabla$) for each method. The results show that the errors decrease with the expected convergence rates, in agreement with the theoretical predictions given in \eqref{eq:apriori1} and \eqref{eq:apriori2}.

For comparison purposes, we also include the errors obtained by solving Problem \eqref{eq:modelproblem} using a well-established polytopal approach, namely the Virtual Element Method \cite{LBe16}, on the same polygonal meshes. Additionally, we report the results obtained with a standard Finite Element Method applied to a set of five Delaunay triangulations \cite{triangle} with mesh sizes $h$ comparable to those of the polygonal grids. The finest mesh of this family is depicted in Figure~\ref{fig:triangle_mesh}.

These results show that our solution converges with the same rate of VEM on the same mesh, and of FEM on a comparable mesh, for orders ranging from $k=1$ to $k=6$, since the solution is regular enough.  
Moreover, it is important to note that, with respect to VEM, the convergence curves indicate that Z-FEM achieves similar accuracy with fewer degrees of freedom. 

To conclude, we observe that the proposed method is more efficient in terms of the number of degrees of freedom than VEM, while maintaining its capability to handle generic polygonal meshes and, therefore, being more flexible than FEM.

\section{Conclusion}
In this paper, we introduce a novel polygonal finite element method called Zipped Finite Element Method (Z-FEM), whose shape functions are given by a linear combination of FEM basis functions, defined on a trivial sub-triangulation of each element of a generic polygonal mesh. 

We present a detailed analysis on the choice of the coefficients of such a combination.
In particular, these coefficients can be efficiently computed by solving a local optimization problem, whose solution can be cheaply obtained by solving a block-diagonal linear system with a single repeated factorizable block. We have demonstrated that, using a particular distribution of degrees of freedom, 
the solution of the optimization problem allows to retrieve all polynomials up to the specified order as a linear combination of the Z-FEM basis functions. We also highlight that, since the explicit form of the involved basis functions is available, this new formulation provides a way to obtain a finite element space that preserves all the theoretical properties of the FEM framework, but on more generic star-shaped polygons. 

Numerical results confirm that our basis functions are able to reproduce polynomials and show optimal performance in terms of standard errors. The proposed method appears to have great potential, combining considerable theoretical simplicity with valuable numerical performance.

\label{sec:conclusion}

\section*{Acknowledgements}

This material is based upon work supported by the Swedish Research Council under grant no. 2021-06594 while the G.T. was in residence at Institut Mittag-Leffler in Djursholm, Sweden during the Fall of 2025.

The author S.B. kindly acknowledges partial financial support provided by European Union through project Next Generation EU, M4C2, PRIN 2022 PNRR project P2022BH5CB\_001 ``Polyhedral Galerkin methods for engineering applications to improve disaster risk forecast and management: stabilization-free operator-preserving methods and optimal stabilization methods'', and by PNRR M4C2 project of CN00000013 National Centre for HPC, Big Data and Quantum Computing (HPC) (CUP: E13C22000990001). The author L.N. kindly acknowledges the financial support provided . The author M.P. kindly acknowledges financial support provided by PEPR/IA (\url{https://www.pepr-ia.fr/}). The author G.T. kindly acknowledges the financial support provided by INdAM-GNCS Project ``Metodi numerici efficienti per problemi accoppiati in sistemi complessi'' (CUP: E53C24001950001). The authors L.N. and G.T. kindly acknowledge financial support provided by project NODES which has received funding from the MUR-M4C2 1.5 of PNRR funded by the European Union - NextGenerationEU (Grant agreement no. ECS00000036) and by the European Union through PRIN project 20227K44ME ``Full and Reduced order modelling of coupled systems: focus on non-matching methods and automatic learning (FaReX)'' (CUP:E53D23005510006). The author F.V. acknowledges the financial support by INdAM-research group GNCS, project title: ``Metodi numerici avanzati per equazioni alle derivate parziali in fenomeni di trasporto e diffusione'' - CUP E53C24001950001.

\bibliographystyle{IEEEtranDOI}
\bibliography{biblio.bib}

\end{document}